\documentclass[11pt]{amsart}
\usepackage{amssymb}
\usepackage{amsmath}
\usepackage{verbatim}

\begin{document}

\newtheorem{theorem}{Theorem}    
\newtheorem{proposition}[theorem]{Proposition}
\newtheorem{conjecture}[theorem]{Conjecture}
\def\theconjecture{\unskip}
\newtheorem{corollary}[theorem]{Corollary}
\newtheorem{lemma}[theorem]{Lemma}
\newtheorem{sublemma}[theorem]{Sublemma}
\newtheorem{observation}[theorem]{Observation}
\theoremstyle{definition}
\newtheorem{definition}{Definition}
\newtheorem{remark}{Remark}
 \def\theremark{\unskip}
\newtheorem{question}[remark]{Question}
\newtheorem{example}{Example}
\def\theexample{\unskip}
\newtheorem{problem}{Problem}

\numberwithin{theorem}{section}
\numberwithin{definition}{section}
\numberwithin{equation}{section}
\numberwithin{remark}{section}

\newcommand{\Norm}[1]{ \left\|  #1 \right\|} 
\newcommand{\norm}[1]{ \|  #1 \|} 
\newcommand{\bignorm}[1]{ \big\|  #1 \big\|}
\newcommand{\abr}[1]{\langle #1 \rangle}
\def\bull{$\bullet$\ } 

\def\jarrow{{\mathbf j}}
 \def\tauvector{{\boldsymbol{\tau}}}
 \def\rhovector{{\boldsymbol{\rho}}}
\def\xivector{\vec{\xi\,}}
\def\rhovector{\vec{\rho}}
\def\tauvector{\vec{\tau}}
\def\Narrow{{\mathbf N}}
\def\Marrow{{\mathbf M}}
\def\karrow{\mathbf k}
\def\rooot{{\mathbf r}}

\def\reals{{\mathbb R}}
\def\torus{{\mathbb T}}
\def\integers{{\mathbb Z}}
\def\complex{{\mathbb C}\/}
\def\naturals{{\mathbb N}\/}
\def\distance{\operatorname{distance}\,}
\def\degree{\operatorname{degree}\,}
\def\kernel{\operatorname{kernel}\,}
\def\dim{\operatorname{dimension}\,}
\def\Span{\operatorname{span}\,}
\def\e{\varepsilon}
\def\p{\partial}
\def\rp{{ ^{-1} }}
\def\Re{\operatorname{Re\,} }
\def\Im{\operatorname{Im\,} }
\def\ov{\overline}
\def\eps{\varepsilon}
\def\lt{L^2}

\def\scriptx{{\mathcal X}}
\def\scriptj{{\mathcal J}}
\def\scriptr{{\mathcal R}}
\def\scripts{{\mathcal S}}
\def\scriptb{{\mathcal B}}
\def\scripta{{\mathcal A}}
\def\scriptk{{\mathcal K}}
\def\scriptd{{\mathcal D}}
\def\scriptp{{\mathcal P}}
\def\scriptl{{\mathcal L}}
\def\scriptv{{\mathcal V}}
\def\scripti{{\mathcal I}}
\def\scriptm{{\mathcal M}}
\def\scripte{{\mathcal E}}
\def\scriptt{{\mathcal T}}
\def\scriptb{{\mathcal B}}
\def\scriptf{{\mathcal F}}
\def\scriptn{{\mathcal N}}
\def\frakg{{\mathfrak g}}
\def\frakG{{\mathfrak G}}
\def\fraks{{\mathfrak S}}
\def\fourierL{{\mathcal F L}}

\author{Michael Christ}
\address{
        Michael Christ\\
        Department of Mathematics\\
        University of California \\
        Berkeley, CA 94720-3840, USA}
\email{mchrist@math.berkeley.edu}
\thanks{The first author was supported by NSF grant DMS-040126}

\author{James Colliander}
\address{James Colliander \\
         Department of Mathematics\\
         University of Toronto\\
       Toronto, ON M5S 2E4, CANADA}
\email{colliander@math.toronto.edu}
\thanks{The second author was supported by NSERC grant RGPIN 250233-03}

\author{Terence Tao}
\address{Terence Tao \\
Department of Mathematics\\
Box 951555 \\
University of California \\
Los Angeles, CA 90095-1555, USA}
\email{tao@math.ucla.edu}
\thanks{The third author was supported by a grant from the MacArthur Foundation.}


\date{December 15, 2006.} 

\title [Nonlinear Schr\"odinger equation]
{A priori bounds and weak solutions
\\
for the nonlinear Schr\"odinger equation
\\
in Sobolev spaces of negative order}

\begin{abstract}
Solutions to the Cauchy problem for the one-dimensional
cubic nonlinear Schr\"odinger equation on the real line 
are studied in Sobolev spaces $H^s$, for $s$ negative but close to $0$.
For smooth solutions 
there is an {\em a priori} upper bound for the $H^s$ norm of
the solution, in terms of the $H^s$ norm of the datum,
for arbitrarily large data, for sufficiently short time.
Weak solutions are constructed for arbitrary initial data in $H^s$.
\end{abstract}

\maketitle

\section{Introduction}
The Cauchy problem for the one-dimensional cubic nonlinear Schr\"odinger equation is
\begin{equation} \label{nls}
\tag{NLS}
\left\{
\begin{aligned}
&iu_t + u_{xx} + \omega |u|^2 u=0
\\
&u(0,x)=u_0(x).
\end{aligned}
\right.
\end{equation}
Here $u=u(t,x)$ with $(t,x)\in[0,T]\times\reals^1$,
and $\omega = \pm 1$.
As is well known, this Cauchy problem is globally wellposed in $H^0$ 
\cite{tsutsumi}.
For all negative $s$ it is illposed in $H^s$, in the sense
that solutions (for smooth initial data) fail to depend {\em uniformly} 
continuously on initial data in the $H^s$ norm
\cite{kpv2},\cite{christcolliandertao}. Moreover,
for $s<-\tfrac12$, there is a stronger form of illposedness:
the solution operator fails even to be continuous
at $0$; there exist smooth solutions with arbitrarily small $H^s$
norms at time $0$, yet arbitrarily large $H^s$ norms
at time $\eps$, for arbitrarily small $\eps>0$.

Our first result, concerning smooth (or more precisely, $H^0$) 
solutions, implies continuity of
the solution map at $u_0=0$
in the $C^0(H^s)$ norm for negative $s$ sufficiently close to $0$, 
in contrast with the strong illposedness for $s<-\tfrac12$.
It asserts an {\em a priori} upper bound for the $H^s$ norm of
an arbitrary smooth solution, 
in terms of the $H^s$ norm of its datum.

\begin{theorem}[A priori bound] \label{thm:aprioribound}
Let $s>-\tfrac1{12}$. Then for all $R<\infty$,
there exist $R'<\infty$ and $T>0$ 
such that for all $u_0\in H^0$ satisfying $\|u_0\|_{H^s}<R$,
the standard solution $u$ of \eqref{nls}
with initial datum $u_0$
satisfies $\max_{t\in[0,T]} \|u(t,\cdot)\|_{H^s}\le R'$.
\end{theorem}

For large $R$, $T$ scales like a certain negative power of $R$.

By the standard solution we mean the unique solution of \eqref{nls}
belonging to the function space $X^{0,b}$ for some $b>\tfrac12$,
or equivalently to $C^0(H^0)\cap L^4([0,T]\times\reals)$.
We know of no reason to believe that the threshold $-\tfrac1{12}$ is optimal; 
slight improvement may be possible via small modifications of our analysis.

Wellposedness of \eqref{nls} has been established by earlier authors
in various function spaces which are wider than $H^0$
\cite{vargasvega},\cite{grunrock},\cite{christ} and scale like
negative order Sobolev spaces, but do not contain $H^s$ for any $s<0$.
We emphasize that those results have a different character than ours;
{\em uniformly} continuous dependence on the initial datum in the 
norm in question is established in those works, whereas it 
certainly fails to hold
\cite{kpv2},\cite{christcolliandertao} in $H^s$ for $s<0$.  

Our second main result asserts the solvability of the Cauchy problem,
in a weak sense, for all initial data in $H^s$ for a range of negative
exponents $s$.
The precise statement involves certain function spaces $Y^{s,b}$, 
which will be specified in Definition~\ref{defn:Ysb}. These are 
variants of the spaces $X^{s,b}$ commonly employed in connection with this equation. 
For any $u\in Y^{s,b}$, $|u|^2 u$ has a natural interpretation as a distribution
(provided that $s>-\tfrac1{12}$ and $b>\tfrac12$),
in the sense that when the space-time Fourier transform
of $|u|^2u$ is written as an integral expression directly in
terms of the space-time Fourier transforms of the factors $u,\bar u,u$,
the resulting integral is absolutely convergent almost everywhere
and defines a tempered locally integrable function; see \eqref{trilinearsbY}.
Thus there is a natural notion of a weak solution in $Y^{s,b}$:
We say that $u\in Y^{s,b}$ 
is a weak solution of \eqref{nls} if the equation holds in the
sense of distributions, when $|u|^2 u$ is interpreted as the inverse
Fourier transform of the function defined by this absolutely convergent
integral.

\begin{theorem}[Existence of weak solutions] \label{thm:weaksolns}
Let $s>-\tfrac1{12}$.
Then there exists $b>\tfrac12$ such that for each $R<\infty$
there exist $R'<\infty$ and $T>0$ 
such that for all $u_0\in H^s$ satisfying $\|u_0\|_{H^s}<R$,
there exists a weak solution $u \in C^0([0,T],H^s)\cap Y^{s,b}$
of \eqref{nls} with initial datum $u_0$ which satisfies $\max_{t \in
  [0,T]} \| u(t) \|_{H^s} \leq R'$.
\end{theorem}


$Y^{s,b}$ embeds continuously in $C^0(H^s)$ for $b>\tfrac12$, so the
$C^0(H^s)$ part of the conclusion is redundant, and is included
only for emphasis.

We do not know whether these weak solutions are unique, let alone whether
there exists any $s<0$ for which the mapping from datum to solution is continuous. 

Our analysis does not rely on the  complete integrability
\cite{ablowitzetal} of \eqref{nls}.
Our arguments would apply, with essentially no changes, to nonintegrable vector-valued
generalizations of the one-dimensional cubic nonlinear Schr\"odinger equation,
provided that those systems obey $H^0$ norm conservation.

We have learned that Theorem~\ref{thm:aprioribound} has been proved independently,
but with a better range of exponents, by H.~Koch and D.~Tataru
\cite{kochtataru}, by arguments
which have much in common with ours.
We are grateful to Justin Holmer for helpful comments.

 \section{Strategy of the analysis}

The strategy is as follows. 
We begin by using the differential equation to (formally, at least) rewrite the increment
$\norm{u(t)}_{H^s}^2-\norm{u_0}_{H^s}^2$
as a multilinear expression in terms of the space-time
Fourier transform of $u$. Certain cancellations
arise, which have no analogues in the corresponding
Fourier expression for $u(t,x)-u_0(x)$.
This leads to an {\it a priori} inequality
of the form 
$\big|\,\norm{u(t)}_{H^s}^2-\norm{u_0}_{H^s}^2\,\big|
\le C\norm{u}_{X^{r,b}}^4,
$
for certain $r,s,b$ with $s<0$ and $r<s$. 
It is this initial step which breaks down if $u$ is replaced by the 
difference of two solutions, preventing us from establishing
any continuity of the map $u_0\mapsto u$.

Thus a bound is required for the $X^{r,b}$ norm, in terms of the $C^0(H^s)$ norm,
but a loss is permitted in the sense that $r$  can be less than $s$.
In \S\ref{section:Ysb}
we introduce certain function spaces $Y^{s,b}$. Their main relevant
properties are:
\begin{enumerate}
\item For $s<0$, $Y^{s,b}$ embeds in $X^{r,b}$, provided that $r<(1+4b)s$.
\item $Y^{s,b}$ embeds in $C^0_t(H^{s-\eps}_x)$ for all\footnote{ 
The definition of $Y^{s,b}$ could be modified so that this would
hold for $\eps=0$, at the expense of small additional
complications in the analysis.}
$\eps>0$, provided that $b>\tfrac12$.
\item If $u,v,w\in Y^{s,b}$ then $u\bar v w\in Y^{s,b-1}$,
under certain restrictions on $s,b$.
\item For solutions of \eqref{nls},
there is an {\it a priori} bound for the $Y^{s,b}$ norm
in terms of the $C^0(H^s)$ norm, of the form
$\norm{u}_{Y^{s,b}}
\le C\norm{u}_{C^0(H^s)}
+ C\norm{u}_{Y^{s,b}}^3$,
valid under certain restrictions on $s,b$. 
\end{enumerate}

Thus one obtains a coupled system of two inequalities
relating $\norm{u}_{C^0(H^s)}$ and $\norm{u}_{Y^{s,b}}$ to
$\norm{u_0}_{H^s}$.
By restricting attention to a short time $T$ and rescaling, one can
reduce matters (for $s>-\tfrac12$) to the case where $u_0$ has small
$H^s$ norm. 
Via a continuity argument, the coupled system then yields a
bound for $\norm{u}_{C^0(H^s)}+\norm{u}_{Y^{s,b}}$ in terms of $\norm{u_0}_{H^s}$.

Weak solutions are obtained as limits of smooth solutions;
an {\it a priori} bound in $H^s$ yields compactness 
in $H^{s-\eps}$ on bounded spatial regions, and it then follows readily 
from the machinery underlying the {\it a priori} bound that
a weak limit of smooth solutions is a weak solution.

An additional argument is needed to place these weak solutions
in $C^0(H^s)$, rather than $C^0(H^{s-\eps})\cap L^\infty(H^s)$.
We refine the machinery by replacing the squared $H^s$ norm 
$\int |\widehat{u}(t,\xi)|^2(1+|\xi|^2)^s\,d\xi$
by 
$\int |\widehat{u}(t,\xi)|^2\varphi(\xi)\,d\xi$
for weight functions $\varphi$ adapted to individual initial data,
so that $\varphi(\xi)\gg (1+|\xi|^2)^s$ for very large $|\xi|$,
and show that control of 
$\int |\widehat{u_0}(\xi)|^2\varphi\,d\xi$ 
extends to control of 
$\int |\widehat{v}(t,\xi)|^2\varphi\,d\xi$ 
for all solutions $v$ of \eqref{nls} with smooth initial data
sufficiently close in $H^s$ norm to $u_0$. This extra control
at high frequencies leads to compactness in $C^0(H^s)$.

\section{Bounding the norm}

In this section we begin to establish an {\em a priori} bound 
for the $C^0(H^s)$ norm of any sufficiently smooth solution of \eqref{nls},
in terms of certain other norms. 
For technical reasons we work with the modified Cauchy problem
\begin{equation} \label{nlsmodified}
\tag{NLS*}
\left\{
\begin{aligned}
&iu_t + u_{xx} + \zeta_0(t)\omega |u|^2 u=0 \\
&u(0,x)=u_0(x)
\end{aligned}
\right.
\end{equation}
where $\zeta_0$ is a smooth real-valued function which is $\equiv 1$
on $[0,T]$, and is supported in $(-2T,2T)$.
Standard proofs of wellposedness in $H^0$ (or in $H^t$ for $t \geq 0$) 
apply to this modified equation.
One advantage is that $u$ can be extended to a solution defined
for all $t\in\reals$.

We will study $\zeta_1(t)u(t,x)$ where $\zeta_1$ is another
real-valued smooth cutoff function supported in $(-2T,2T)$
which satisfies $\zeta_1\zeta_0\equiv\zeta_0$.
Because the equation is simply
the linear Schr\"odinger equation outside the support of
$\zeta_0$, a $C^0(H^s)$ bound holds for $\zeta(t)u(t,x)$
for one real-valued cutoff function in $\zeta\in C^\infty_0(-2T,2T)$
satisfying $\zeta\zeta_0\equiv\zeta_0$
if and only if such a bound holds for every such function.

Recall \cite{bourgain} the function space
$X^{s,b}$, which is defined to be the set of all space-time distributions $u$ whose spacetime Fourier transform $\widehat{u}$ is such that
\begin{equation*}
\norm{u}_{X^{s,b}}^2
:= \iint_{\reals^2} |\widehat{u}(\xi,\tau)|^2
\abr{\tau-\xi^2}^{2b}\abr{\xi}^{2s}\,d\xi\,d\tau
<\infty
\end{equation*}
where of course $\abr{x} := (1 + |x|^2)^{1/2}$.

One of the two principal inequalities underlying our theorems is
as follows. The second is formulated in Proposition~\ref{prop:Ysbbound}.
\begin{proposition} \label{prop:PhifromX}
Let $T_0<\infty$, $T\in[0,T_0]$, $s\in (-\tfrac12,0)$, $b\in (\tfrac12,1)$.
There exists $C<\infty$ such that
for any sufficiently smooth solution\footnote{For instance, $u_0 \in H^{10}$ would suffice.} $u$ of \eqref{nlsmodified}
with initial datum $u_0$,
\begin{equation}
\Big| \Norm{u}_{C^0([-2T,2T],H^s)}^2 - \Norm{u_0}_{H^s}^2 \Big|
\le C\Norm{\zeta_1 u}_{X^{r,b}}^4
\end{equation}
provided that
\begin{equation} \label{PhifromXrrestriction}
r>-\tfrac14\ \text{and}\  b>\tfrac12.
\end{equation}
\end{proposition}
For $s>-\tfrac14$, the right-hand side
involves a norm which is weaker, in terms of the number of spatial
derivatives involved, than the $C^0(H^s)$ norm.
The proof of this result is begun below and 
completed in \S\ref{section:conclusionofPhifromX},
using some of the inequalities established in \S\ref{section:strichartz}.

We will work with both spatial Fourier coefficients
\begin{equation}
\widehat{u}(t,\xi) := \int_{\reals} e^{-ix\xi}u(t,x)\,dx
\end{equation}
and space-time Fourier coefficients
\begin{equation}
\widehat{u}(\xi,\tau) := \int_{\reals^2} e^{-ix\xi}e^{-it\tau}u(t,x)\,dx\,dt;
\end{equation}
it will be clear from context and from the names of the variables which
of these two is meant in any particular instance.
The differential equation \eqref{nlsmodified}
is expressed in terms of spatial Fourier coefficients as
\begin{equation}
\frac{d}{dt} \widehat{u}(t,\xi)
= -i\xi^2 \widehat{u}(t,\xi) + i\omega\zeta_0(t)\int_{\xi_1-\xi_2+\xi_3=\xi} 
\widehat{u}(t,\xi_1)
\overline{\widehat{u}(t,\xi_2)}
\widehat{u}(t,\xi_3)
\,d\lambda_\xi
\end{equation}
where $\lambda_\xi$ is appropriately normalized Lebesgue
measure on $\{(\xi_1,\xi_2,\xi_3)\in\reals^3: \xi_1-\xi_2+\xi_3=\xi\}$.

Consider any sufficiently regular solution $u$ of \eqref{nlsmodified}.
Let $\varphi:\reals\to[0,\infty)$ and define the modified mass
\begin{equation}
\Phi_\varphi(t) = 
\Phi_\varphi(t,u) := \int_\reals |\widehat{u}(t,\xi)|^2\varphi(\xi)\,d\xi.
\end{equation}
We will be primarily interested in $\varphi(\xi) = \abr{\xi}^{2s}$,
but more general weights will be needed to establish the full
conclusion of Theorem~\ref{thm:weaksolns}.

A short calculation shows that we have the ``almost conservation law''
$$\frac{d\Phi}{dt} = \Re(c\omega \scripti)$$
for $\Phi$, where $c$ is an absolute constant, $\scripti$ is the multilinear integral
\begin{equation} \label{Phiprime}
\scripti(t)
=\scripti_\varphi(u,t) := 
\zeta_0(t)\int_\Xi 
\widehat{u}(t,\xi_1)
\overline{\widehat{u}(t,\xi_2)}
\widehat{u}(t,\xi_3)
\overline{\widehat{u}(t,\xi_4)}
\,\psi(\xivector)
\,d\lambda(\xivector),
\end{equation}
$\xivector = (\xi_1,\cdots,\xi_4)\in\reals^4$ is a multi-frequency, 
$\Xi \subset \reals^4$ is the hyperplane
\begin{equation}
\Xi := \{\xivector: \xi_1-\xi_2+\xi_3-\xi_4=0\},
\end{equation}
$\lambda$ is appropriately normalized Lebesgue measure on $\Xi$,
and
\begin{equation}
\psi(\xivector)
:= \varphi(\xi_1)
-\varphi(\xi_2)
+\varphi(\xi_3)
-\varphi(\xi_4).
\end{equation}
Thus\footnote{As usual, we use $X \lesssim Y$ to denote an estimate of the form $X \leq CY$ for some constant $C$, depending only on the exponents $r$, $s$ and $b$ which will appear later in this paper.} 
$|\Phi(t)-\Phi(0)| \lesssim 
|\int_0^t \scripti(r)\,dr|$.

Introduce also 
\begin{equation}
\sigma(\xi_1,\cdots,\xi_4) := 
\xi_1^2-\xi_2^2+\xi_3^2-\xi_4^2.
\end{equation}
$\sigma$ has the useful alternative expressions 
\begin{equation}
\sigma(\xivector) = 2(\xi_1-\xi_2)(\xi_3-\xi_2) 
=2 (\xi_1-\xi_4)(\xi_3-\xi_4)
=2 (\xi_1-\xi_2)(\xi_1-\xi_4)
\ \ \forall\,\xivector\in \Xi.
\end{equation}

We have the following basic cancellation bound (cf. \cite{ckstt}):

\begin{lemma}[Double mean value theorem]\label{lemma:psibound}
Let $\xivector=(\xi_1,\cdots,\xi_4)\in\Xi\subset\reals^4$.
If $\varphi\in C^2$ and all $\xi_j$ belong to a common interval $I$ then
$|\psi(\xivector)| \le |\sigma(\xivector)| \max_{y\in I} |\varphi''(y)|$.
\end{lemma}

\begin{proof}
$\varphi(\xi_2)-\varphi(\xi_1)
= (\xi_2-\xi_1)\int_0^1 \varphi'(\xi_1+t(\xi_2-\xi_1))\,dt$.
Writing the corresponding expression for $\varphi(\xi_4)-\varphi(\xi_3)$,
and noting that $(\xi_2-\xi_1) = -(\xi_4-\xi_3)$ 
since $\xivector\in\Xi$, gives
\begin{align*}
\psi(\xi)
&= (\xi_2-\xi_1)\int_0^1 [\varphi'(\xi_1+t(\xi_2-\xi_1))
- \varphi'(\xi_4+t(\xi_3-\xi_4))]\,dt
\\
&=
(\xi_2-\xi_1)(\xi_1-\xi_4)
\iint_{[0,1]^2} \varphi''(\xi_1+t(\xi_2-\xi_1)+s(\xi_4-\xi_1)\,ds\,dt.
\end{align*}
\end{proof}

In order to control the contribution made by the region not close to the
diagonal,
express each factor $\widehat{u}(t,\xi)$ in the integral
as the inverse Fourier transform of its Fourier transform with respect to $t$,
to obtain for all $t\in[-2T,2T]$
\begin{equation} \label{Phiincrementbound}
|\int_0^t \scripti_\varphi(u,r)\,dr|
\le C
\int_{\Xi}
\int_{\reals^4}
\prod_{j=1}^4 |\widehat{u}(\xi_j,\tau_j)|
\abr{\tau_1-\tau_2+\tau_3-\tau_4}^{-1}
|\psi(\xivector)|
\,d\tauvector
\,d\lambda(\xivector)
\end{equation}
where $C$ depends on $T$ and
$\tauvector=(\tau_1,\cdots,\tau_4)$.
The notation $\widehat{u}$ denotes here the Fourier transform
with respect to both spatial and temporal variables.

Write
\begin{equation} \label{substitutegsubn}
|\widehat{u}(\xi_j,\tau_j)|
=: \abr{\xi_j}^{-r}\abr{\tau_j-\xi_j^2}^{-b}
g_j(\xi_j,\tau_j). 
\end{equation}
Then $\norm{g_j}_{L^2(\reals)}=\norm{u}_{X^{r,b}}$.
The right-hand side of \eqref{Phiincrementbound} becomes
\begin{equation} \label{needStrichartz}
\int_{\reals^4}
\int_\Xi
\prod_{n=1}^4 \Big(
g_n(\xi_n,\tau_n)
\abr{\xi_n}^{-r}\abr{\tau_n-\xi_n^2}^{-b}
\Big)
|\psi(\xivector)|
\,d\lambda(\xivector) 
\abr{\tau_1-\tau_2+\tau_3-\tau_4}^{-1}
\,d\tauvector.
\end{equation}
In the \S\ref{section:conclusionofPhifromX}
we will complete the proof of Proposition~\ref{prop:PhifromX}
by showing that
for $\varphi(\xi) = \abr{\xi}^{2s}$, the integral \eqref{needStrichartz} 
is majorized by
$C\prod_{n=1}^4\norm{g_n}_{\lt(\reals^2)}$
provided that $s,r,b$ satisfy the hypotheses of the proposition.

\section{Trilinear inequalities of Strichartz type}
\label{section:strichartz}

A prototypical inequality of Strichartz type
says that for $h\in\lt(\reals)$, 
the solution $u$ of the linear Schr\"odinger equation with
initial datum $h$ belongs to $L^6(\reals^2)$. 
Therefore any three such solutions satisfy $u_1\bar u_2 u_3\in\lt$.
Rewritten on the Fourier side by means of the Plancherel identity,
this becomes 
\begin{equation} \label{trivialStrichartz}
\big|\int f(\xi_1-\xi_2+\xi_3,\xi_1^2-\xi_2^2+\xi_3^2)
\prod_{n=1}^3 g_n(\xi_n) \,d\xi_n\big|
\lesssim \prod_{n=1}^3 \norm{g_n}_{L^2(\reals)}\norm{f}_{\lt(\reals^2)}.
\end{equation}
One version of the bilinear Strichartz inequality,
expressed directly in terms of Fourier variables,
states that for any subset $E\subset\reals^2$,
\begin{multline} \label{bilinearStrichartz}
\Big|\int_{\reals^2} f(\xi_1\pm \xi_2,\xi_1^2\pm\xi_2^2)
h_1(\xi_1)h_2(\xi_2)
\chi_E(\xi_1,\xi_2)
\,d\xi_1\,d\xi_2 \Big|
\\
\lesssim \big(\min_{(\xi_1,\xi_2)\in E}|\xi_1-\xi_2|\big)^{-1/2}
\norm{f}_{L^2(\reals^2)}
\norm{h_1}_{L^2(\reals^1)}
\norm{h_2}_{L^2(\reals^1)},
\end{multline}
where the two $\pm$ signs are either both $+$, or both $-$;
this represents the pairing of $f$ with a bilinear operator applied
to $h_1,h_2$.
This is due to Carleson and Sj\"olin \cite{carlesonsjolin},
and is a direct consequence of Cauchy-Schwarz via the substitution
$(\xi_1,\xi_2)\mapsto  (\xi_1\pm\xi_2,\xi_1^2\pm\xi_2^2)$.
Its advantage, in practice, is that it provides a superior bound
when $|\xi_1-\xi_2|$ is large. 

In this section we establish certain versions of the trilinear
inequality \eqref{trivialStrichartz} which incorporate improvements similar
to the factor $|\xi_1-\xi_2|^{-1/2}$ in \eqref{bilinearStrichartz}.
These arise naturally in the analysis of the Fourier
transform of a threefold product $u\bar v w$ of 
functions in spaces $X^{r,b}$ or $Y^{s,b}$.

\subsection{Statements of inequalities}
\begin{proposition} \label{prop:strichartzA}
Consider
\begin{equation} \label{strichartzAintegral}
\int_{\xivector\in S\subset \Xi}
\int_{\tau_1-\tau_2+\tau_3-\tau_4=0}
\prod_{n=1}^4  g_n(\xi_n,\tau_n)\abr{\tau_n-\xi_n^2}^{-\beta_n}
\chi_E(L(\xivector))
\,d\lambda(\tauvector)
\,d\lambda(\xivector)
\end{equation}
where each $g_n\ge 0$, 
$i,j\in\{1,2,3,4\}$ are distinct,
$E\subset\reals^1$ is any measurable set,
and $L:\reals^4\to\reals$ is a linear transformation.
Suppose that
\newline
\bull $\beta_n>\tfrac12$ for all but at most one index $n$,
and $\beta_n>0$ for all $n$.
\newline
\bull $i,j$ have opposite parity.
\newline
\bull $L$ 
belongs neither to the linear span of
$\{\xi_i,\xi_j, \xi_1-\xi_2+\xi_3-\xi_4\}$,
nor to the linear span of
$\{\xi_k,\xi_l, \xi_1-\xi_2+\xi_3-\xi_4\}$,
where $\{i,j,k,l\}=\{1,2,3,4\}$.

Then there exists $C<\infty$ depending on $L$ such that
\eqref{strichartzAintegral} is majorized by
\begin{equation} \label{strichartzAbound}
C
|E|^{1/2}
\max_{\xivector \in S}
\big(
|\xi_i-\xi_j|^{-1/2}
\abr{\sigma(\xivector)}^{-\beta}
\big)
\prod_{n=1}^4\norm{g_n}_{L^2(\reals^2)}
\end{equation}
where $\beta=\min_n\beta_n$. 
\end{proposition}

In our application, $L$ will take the form
$L(\xivector) = \xi_\mu-\xi_\nu$ for some $\mu\ne\nu$.
If $\{\mu,\nu\}$ equals neither
$\{i,j\}$ nor $\{1,2,3,4\}\setminus\{i,j\}$
then $L$ satisfies the hypothesis.

A variant of this inequality applies to other
linear transformations $L$:
\begin{proposition} \label{prop:strichartzC}
Consider \eqref{strichartzAintegral}
with $L(\xivector)=\xi_k$ for some $k\notin\{i,j\}$.
Suppose again that
$\beta_n>\tfrac12$ for all but at most one index $n$,
and $\beta_n>0$ for all $n$.
Suppose that $|E|\lesssim
\min_{\xivector\in S} |\xi_i-\xi_j|$.
Let $\beta :=\min_n\beta_n$.
Then there exists $C<\infty$ such that
\eqref{strichartzAintegral} is majorized by
\begin{equation} \label{strichartzCbound}
\lesssim |E|^{1/4}
\max_{\xivector \in S} 
\big(|\xi_i-\xi_j|^{-1/4} \abr{\sigma(\xivector)}^{-\beta}\big)
\prod_{n=1}^4\norm{g_n}_{L^2(\reals^2)}.
\end{equation}
\end{proposition}

\begin{remark}[Trilinear Knapp example] 
\eqref{strichartzCbound} is (in practice) weaker than \eqref{strichartzAbound},
because $|E|/\min_S|\xi_i-\xi_j|$ is raised
only to the power $\tfrac14$, rather than
$\tfrac12$ as in Proposition~\ref{prop:strichartzA}.
The exponent $\tfrac14$ in \eqref{strichartzCbound} 
is however optimal.
We will show this for  
\begin{equation} \label{simplifiedexample}
\int_{\reals^3} 
G(\xi_1+\xi_3-\xi_4,\xi_1^2+\xi_3^2-\xi_4^2)
\prod_{n\ne 2} h_n(\xi_n)\
\chi_{|\xi_4|\le 1}\chi_S(\xivector) 
\,d\xi_1\,d\xi_3\,d\xi_4,
\end{equation}
which is a simplified version of an expression 
arising in the proof of Proposition~\ref{prop:strichartzC}
(see the case $\nu=2$); the example can be 
adapted to precisely the situation arising there. 

Suppose that $\min_S|\xi_1-\xi_2|\gtrsim 1$. Then
\eqref{simplifiedexample} is bounded by 
$\lesssim \min_{S}|\xi_1-\xi_2|^{-1/4}\norm{G}_{\lt}\prod_{n\ne 2}\norm{h_n}_{\lt}$,
as will be shown in the proof of \eqref{strichartzCbound} below; 
we claim now that the exponent $\tfrac14$ cannot be improved.
To see this,
define $h_4$ to be the characteristic function of $[0,1]$,
$h_1$ to be the characteristic function of the interval $[N,N+\Delta]$,
and
$h_3$ to be the characteristic function of $[N+N^{1/2},N+N^{1/2}+\Delta]$.
Let $A$ be a sufficiently large positive constant, and define
$G(x,y)$ to be the characteristic function of the 
set of all $(x,y)\in\reals^2$ satisfying
$|x-2N|\le AN^{1/2}$ and
$|y+2N^2+2N^{3/2}-2Nx|\le AN$.

A short calculation shows that if $A$ is chosen to be sufficiently large
then 
$G(\xi_1+\xi_3-\xi_4,\xi_1^2+\xi_3^2-\xi_4^2)\equiv 1$ whenever
$h_n(\xi_n)\ne 0$ for all $n\in\{1,3,4\}$.
Therefore the integral is simply
$\prod_{n\ne 2}\int_{\reals^1} h_n\sim N^{1/2}\cdot N^{1/2}\cdot 1=N$.
On the other hand, 
$\norm{G}_{\lt}\sim  N^{3/4}$,
while
$\norm{h_n}_{\lt}\sim N^{1/4}$
for $n=1,3$
and $\sim 1$ for $n=4$.
Thus the product of the four $\lt$ norms has order of magnitude
$N^{5/4}$, and consequently
the ratio of \eqref{simplifiedexample} to the product of norms
has order of magnitude $N/N^{5/4} = N^{-1/4}$. 
Since $\xi_2-\xi_1=\xi_3-\xi_4$ has order of magnitude $N$,
this is the ratio claimed. \qed
\end{remark}

The analysis of \eqref{needStrichartz}
is a bit more complicated because the relation
$\tau_1-\tau_2+\tau_3-\tau_4=0$ is 
replaced by the slowly decaying factor
$\abr{\tau_1-\tau_2+\tau_3-\tau_4}^{-1}$.
It requires a third variant:
\begin{proposition} \label{prop:strichartzB}
Consider
\begin{equation} \label{strichartzBintegral}
\int_{\xivector\in S\subset \Xi}
\int_{\tauvector\in\reals^4}
\abr{\tau_1-\tau_2+\tau_3-\tau_4}^{-1}\,
\prod_{n=1}^4  g_n(\xi_n,\tau_n)\abr{\tau_n-\xi_n^2}^{-\beta_n}
\phi(\xivector)
\chi_E(L(\xivector))
\,d\lambda(\xivector)\,d\tauvector
\end{equation}
where $g_n\ge 0$, $\phi\ge 0$, 
and $\beta_n>\tfrac12$ for all $n$.
Let $i\ne j\in\{1,2,3,4\}$ and 
let $L:\reals^4\to\reals$ be a linear functional satisfying 
the hypotheses of Proposition~\ref{prop:strichartzA}.
Then \eqref{strichartzBintegral} is majorized by
\begin{equation} \label{strichartzBbound}
\lesssim
\prod_{n=1}^4\norm{g_n}_{L^2(\reals^2)}
|E|^{1/2}
\max_{\xivector\in S}
\Big[
\phi(\xivector) 
|\xivector|^{1/2}
\abr{\sigma(\xivector)}^{-1} 
+ 
\phi(\xivector) 
|\xi_i-\xi_j|^{-1/2}
\abr{\sigma(\xivector)}^{-\beta}
\Big]
\end{equation}
where $\beta := \min_n\beta_n$.
\end{proposition}

While \eqref{strichartzBintegral} is formally similar to \eqref{strichartzAintegral}, 
a significant contribution to \eqref{strichartzBintegral}
can arise from a region in Fourier space which has no analogue
in \eqref{strichartzAintegral}. 
This region contributes an additional term in \eqref{strichartzBbound}.
In our application, $\phi$ will be $|\psi|$.

\subsection{Proofs of inequalities}
The essence of 
Propositions~\ref{prop:strichartzA}, \ref{prop:strichartzC},
and \ref{prop:strichartzB}
lies in the following two simpler inequalities.
\begin{lemma} \label{lemma:boileddownstrichartz}
Let $i,j,k$ be the three elements of $\{1,2,3\}$, written in any order.
Let $\ell:\reals^3\mapsto\reals^1$
satisfy $\partial\ell/\partial\xi_k\ne 0$. Then
for any nonnegative measurable functions $G,g_n$ of two and one real variables
respectively,
any measurable sets $E\subset\reals^1$ and $S\subset\Xi$,
the quantity
\begin{equation}
\int_{\reals^3} \prod_{n=1}^3 g_n(\xi_n) G(\xi_1-\xi_2+\xi_3,
\xi_1^2-\xi_2^2+\xi_3^2)
\chi_S(\xivector)\chi_E(\ell(\xivector)) 
\,d\xi_1\,d\xi_2\,d\xi_3
\end{equation}
is majorized by
\begin{equation}
\lesssim \norm{G}_{\lt}\prod_{n=1}^3 \norm{g_n}_{\lt} |E|^{1/2} 
(\min_{\xivector\in S}|\xi_i-\xi_j|)^{-1/2}
\end{equation}
where the implied constant depends on $\ell$.
\end{lemma}

\begin{proof}
Consider the case where $\{i,j\}=\{1,2\}$.
Apply Cauchy-Schwarz to majorize by
\begin{multline}
\Big(\int_{\xivector\in S} G^2(\xi_1-\xi_2+\xi_3,
\xi_1^2-\xi_2^2+\xi_3^2) g_3^2(\xi_3) 
\,d(\xi_1,\xi_2,\xi_3)
\Big)^{1/2}
\\
\Big(\int_{\reals^3} g_1^2(\xi_1)g_2^2(\xi_2)\chi_E(\ell(\xivector)) 
\,d(\xi_1,\xi_2,\xi_3)
\Big)^{1/2}.
\end{multline}

The left-hand factor is majorized by
$\lesssim \norm{G}_{\lt} \norm{g_3}_{\lt} (\min_S|\xi_1-\xi_2|)^{-1/2}$;
this is seen by first fixing $\xi_3$ and integrating with
respect to $(\xi_1,\xi_2)$, making the change of variables
$(\xi_1,\xi_2)\mapsto (\xi_1-\xi_2,\xi_1^2-\xi_2^2)$.

To analyze the right-hand factor, first integrate with respect to $\xi_3$,
obtaining a bound of 
\begin{equation}
\lesssim |E|^{1/2} (\int g_1^2(\xi_1)g_2^2(\xi_2)\,d\xi_1\,d\xi_2)^{1/2}
\end{equation}
since
$\partial \ell/\partial\xi_3\ne 0$. 
Then integrate with respect to $(\xi_1,\xi_2)$. 
Multiplying these bounds for the two factors yields
$\lesssim \norm{G}_{\lt} \prod_{n=1}^3 \norm{g_n}_{\lt} 
|E|^{1/2} (\min_S|\xi_1-\xi_2|)^{-1/2}$.

The same reasoning applies for other $\{i,j\}$;
in all cases $|\xi_i-\xi_j|$ arises, rather than $|\xi_i+\xi_j|$.
\end{proof}

\begin{lemma} \label{lemma:routine}
For $m=1,2$ let $L_m:\reals^4\to\reals$
be linear functionals such that  
$\{\xi_1-\xi_2+\xi_3-\xi_4,L_1,L_2\}$
is linearly independent.
Then all nonnegative measurable functions
$g_n\in\lt(\reals^1)$ and all measurable sets
$E_m\subset\reals^1$,
\begin{equation}
\int_\Xi
\prod_{n=1}^4 g_n(\xi_n)
\prod_{m=1}^2 
\chi_{E_m}(L_m(\xivector))
\,d\lambda(\xivector)
\lesssim
\prod_{m=1}^2|E_m|^{1/2}
\prod_{n=1}^4 \norm{g_n}_{\lt}.
\end{equation}
\end{lemma}

\begin{proof}
Consider the multlinear form
$T(g_1,\cdots,g_6)
:=\int_\Xi
\prod_{n=1}^4 g_n(\xi_n)
g_5(L_1(\xivector))
g_6(L_2(\xivector))
\,d\lambda(\xivector)$.
By Cauchy-Schwarz,
\begin{equation*}
|T(g_1,\cdots,g_6)|
\le 
\Big( 
\int_\Xi
\prod_{n=1}^3 |g_n(\xi_n)|^2
\,d\lambda(\xivector)
\Big)^{1/2}
\Big( 
\int_\Xi
|g_4(\xi_4)|^2
|g_5(L_1(\xivector))|^2
|g_6(L_2(\xivector))|^2
\,d\lambda(\xivector)
\Big)^{1/2}.
\end{equation*}
The first factor is a constant multiple of
$\prod_{j=1}^3\norm{g_j}_2$.
The assumption that
$\{\xi_1-\xi_2+\xi_3-\xi_4,L_1,L_2\}$
is linearly independent implies that
the second factor is a constant multiple of
$\prod_{j=4}^6\norm{g_j}_2$.
\end{proof}

\begin{proof}[Proof of Proposition~\ref{prop:strichartzA}]
Now consider the quantity given in the statement of the proposition.
Introduce $\rho_n = \tau_n-\xi_n^2$.
Since $\rho_1-\rho_2+\rho_3-\rho_4
= \xi_1^2
- \xi_2^2
+ \xi_3^2
- \xi_4^2
= \sigma(\xivector)$,
we have $\abr{\rho_n}\gtrsim \abr{\sigma(\xivector)}$ for some $n$.
Partition the region of integration into four subregions,
according to the index $n$ for which $|\rho_n|$ is largest.
By symmetry, it suffices to prove the stated bound for one 
of these subregions. 
Let $\nu\in\{1,2,3,4\}$ be arbitrary, and consider the subregion
consisting of all $\xivector$ satisfying
$|\rho_\nu(\xivector)|=\max_n|\rho_n(\xivector)|$.

Suppose first that $\nu\notin\{i,j\}$.
Then choose $\mu$ so that $\{1,2,3,4\}=\{i,j,\mu,\nu\}$.
The contribution of the subregion under examination is 
\begin{multline} \label{eq:subregion}
\lesssim 
\int_{(\rho_1,\rho_2,\rho_3)\in\reals^3}
\Big(
\int_{\reals^3} g_\nu(\pm \xi_\mu\pm\xi_i\pm\xi_j,
\pm \xi_\mu^2\pm\xi_i^2\pm\xi_j^2 \pm\rho_\mu\pm\rho_i\pm\rho_j)
\prod_{n\ne\nu} h_n(\xi_n,\rho_n)
\\
\abr{\rho_\nu}^{-\beta_\nu}
\chi_S(\xivector)\chi_E(L(\xivector))
\,d\xi_i\,d\xi_j\,d\xi_\mu
\Big)
\prod_{n\ne\nu}\abr{\rho_n}^{-\beta_n}
\,d(\rho_i,\rho_j,\rho_\mu)
\end{multline}
where the $\pm$ sign preceding $\xi_n^2$ agrees with the sign
preceding $\xi_n$ for each $n\in\{\mu,i,j\}$,
and where the outer integral extends only over those $\rho_n$
satisfying $|\rho_n|\lesssim \Lambda$.
Here 
$h_n(\xi_n,\rho_n)=g_n(\xi_n,\tau_n) = g_n(\xi_n,\rho_n+\xi_n^2)$,
and consequently
$\norm{h_n}_{\lt}=\norm{g_n}_{\lt}$.

Fix $(\rho_i,\rho_j,\rho_\mu)$.
The linear transformation $\Xi\in\xivector\mapsto 
(\xi_i,\xi_j,\xi_\mu)\in\reals^3$ is invertible,
so there is a unique linear functional $\tilde L(\xivector):\reals^3\mapsto\reals$
satisfying $\tilde L(\xi_i,\xi_j,\xi_\mu)=L(\xivector)$.
The hypothesis on $L$ ensures that $\partial\tilde L/\partial\xi_\mu\ne 0$.
The inner integral thus takes the form discussed in 
Lemma~\ref{lemma:boileddownstrichartz}, 
and is consequently majorized by
\begin{equation}
\lesssim \norm{g_\nu}_{\lt(\reals^2)}
\prod_{n\ne\nu} \norm{g_n(\cdot,\rho_n)}_{\lt(\reals^1)} 
|E|^{1/2} 
\sup_S(|\xi_i-\xi_j|^{-1/2}\abr{\sigma(\xivector)}^{-\beta_\nu})
\end{equation}
since $\abr{\rho_\nu}\gtrsim \abr{\sigma(\xivector)}$
at every point of the region of integration.

It remains to bound
$\prod_{n\ne\nu}
\int_{|\rho_n|\lesssim\Lambda} \norm{g_n(\cdot,\rho_n)}_{\lt(\reals^1)}
\abr{\rho_n}^{-\beta_n}\,d\rho_n$.
If $\beta_n>\tfrac12$ then 
\begin{equation}
\int_{\reals} \norm{g_n(\cdot,\rho_n)}_{\lt(\reals^1)}
\abr{\rho_n}^{-\beta_n}\,d\rho_n
\lesssim \norm{g_n}_{\lt(\reals^2)}
\end{equation}
by Cauchy-Schwarz.
If $\beta_n>\tfrac12$ for all $n\ne\nu$ then the desired bound is obtained.

Otherwise there remains exactly
one index $m\ne\nu$ such that $\beta_m\le\tfrac12$.
Then $\beta=\beta_m$, and $\beta_\nu\ge\beta_m$.
Since $\abr{\rho_\nu}\gtrsim \abr{\rho_m}, \abr{\sigma(\xivector)}$
throughout the region of integration,
one has
\begin{equation*}
\abr{\rho_\nu}^{-\beta_\nu} 
\lesssim \abr{\sigma(\xivector)}^{-\beta_m}
\abr{\rho_m}^{\beta_m-\beta_\nu}
= 
\abr{\sigma(\xivector)}^{-\beta}
\abr{\rho_m}^{\beta_m-\beta_\nu}.
\end{equation*}
This factor of $\abr{\rho_m}^{\beta_m-\beta_\nu}$,
multiplied by the factor of $\abr{\rho_m}^{-\beta_m}$
already present in the integral, becomes $\abr{\rho_m}^{-\beta_\nu}$.
Since $\beta_\nu>\tfrac12$, the analysis can be completed as above.
This concludes the analysis, in the case where $\nu\notin\{i,j\}$.

Suppose finally that $\nu\in\{i,j\}$; by symmetry, we may suppose
that $\nu=i$. Writing $\{1,2,3,4\}=\{i,j,k,l\}$,
the equation for $\Xi$, together with the hypothesis that $i,j$
have opposite parity, imply that
$|\xi_i-\xi_j|\equiv|\xi_k-\xi_l|$ for all $\xivector\in\Sigma$,
whence $\min_S|\xi_i-\xi_j|=\min_S|\xi_k-\xi_l|$.
Therefore $\{i,j\}$ can be interchanged with 
$\{k,l\}$. Define $\mu$ so that $\{1,2,3,4\}=\{\mu,\nu,k,l\}$.
The hypothesis on $L$ is explicitly formulated so as to be unaffected
under this symmetry.
Therefore the above reasoning applies, and again yields the stated bound.
\end{proof}

\begin{proof}[Proof of Proposition~\ref{prop:strichartzC}]
\eqref{needStrichartz} is
invariant under the permutations
$(1,2,3,4)\mapsto (2,1,4,3)$,
$(1,2,3,4)\mapsto (3,2,1,4)$,
$(1,2,3,4)\mapsto (1,4,3,2)$,
and consequently also
$(1,2,3,4)\mapsto (3,4,1,2)$
of the indices.
Therefore it is no loss of generality to assume that $i=1$, $j=2$, and $k=4$.

We follow the proof of Proposition~\ref{prop:strichartzA}.
In the case when $\nu\notin\{1,2\}$, because $L(\xivector)=\xi_4$ does not
belong to the span of the three linear transformations $\xi_1$,
$\xi_2$, and $\xi_1-\xi_2+\xi_3-\xi_4$, that proof 
applies without alteration and yields the upper bound
\eqref{strichartzAbound}. Since $|E|/\min_S|\xi_1-\xi_2|\lesssim 1$
by hypothesis,
\eqref{strichartzAbound} is majorized by a constant multiple
of the desired bound \eqref{strichartzCbound}.

Consider next the case where $\nu=2$.
Then because $\xi_1-\xi_2\equiv -(\xi_3-\xi_4)$,
Lemma~\ref{lemma:boileddownstrichartz}
can be applied with the roles of the indices $2,4$ interchanged
to obtain a bound
\begin{equation}  \label{firstCbound}
\lesssim |E|^{1/2}
\max_{\xivector \in S} 
\big(|\xi_1-\xi_3|^{-1/2}
\abr{\sigma(\xivector)}^{-\beta}\big)
\prod_{n=1}^4\norm{g_n}_{L^2(\reals^2)}.
\end{equation}
Another bound is also available. Apply Proposition~\ref{prop:strichartzA}
with $L$ replaced by $\tilde L(\xivector) = \xi_1-\xi_3$;
$\tilde L$ does not belong to the span of 
$\{\xi_1,\xi_2,\xi_1-\xi_2+ \xi_3-\xi_4\}$, 
nor to the span of
$\{\xi_3,\xi_4,\xi_1-\xi_2+ \xi_3-\xi_4\}$, 
so the hypotheses are satisfied. This yields an
alternative bound
\begin{equation} \label{secondCbound}
\lesssim \max_{\xivector\in S} |\xi_1-\xi_3|^{1/2}
\max_{\xivector \in S} 
\big(|\xi_1-\xi_2|)^{-1/2}
\abr{\sigma(\xivector)}^{-\beta}\big)
\prod_{n=1}^4\norm{g_n}_{L^2(\reals^2)}.
\end{equation}

If $\max_S |\xi_1-\xi_3|$ is comparable to $\min_S |\xi_1-\xi_3|$,
then taking the geometric mean of these two upper bounds yields
the desired bound \eqref{strichartzCbound}.
Decomposing $S$ into subsets $S_\kappa$
in which $|\xi_1-\xi_3|$ is comparable
to $2^\kappa$ for arbitrary $\kappa\in\integers$,
invoking whichever of \eqref{firstCbound},
\eqref{secondCbound} is more favorable for each $\kappa$,
and summing over $\kappa$ yields the same bound in the general case.

Finally, when $\nu=1$,
apply Lemma~\ref{lemma:boileddownstrichartz}
with the roles of the indices $1,4$ interchanged,
and repeat the above discussion for the case $\nu=2$,
replacing $\xi_1-\xi_3$ by $\xi_2-\xi_3$ throughout.
The reasoning is otherwise unchanged.
\end{proof}

\begin{proof}[Proof of Proposition~\ref{prop:strichartzB}]
Substitute $\tau_n=\rho_n+\xi_n^2$ and $g_n(\xi_n,\tau_n)
= \tilde g_n(\xi_n,\rho_n)$  for all $n\in\{1,2,3,4\}$
to transform the integral into
\begin{equation} \label{eq:strichartzB1}
\int_{\reals^4}\int_{S\subset \Xi}
\abr{\rho_1-\rho_2+\rho_3-\rho_4+\sigma(\xivector)}^{-1}
\prod_{n=1}^4 \tilde g_n(\xi_n,\rho_n)\abr{\rho_n}^{-\beta_n}
\phi(\xivector)
\chi_E(L(\xivector))
\,d\lambda(\xivector)\,d\rhovector,
\end{equation}
where $\tilde g_n$ has the same $\lt$ norm as $g_n$. 

Begin with
the region where $|\rho_n|\le \tfrac18\abr{\sigma(\xivector)}$
for all $n\in\{1,2,3,4\}$, which 
has no counterpart in Proposition~\ref{prop:strichartzA}.
Its contribution is comparable to
\begin{equation}
\int_{\reals^4}
\int_{S\subset \Xi}
\prod_{n=1}^4 \tilde g_n(\xi_n,\rho_n)\abr{\rho_n}^{-\beta_n}
\abr{\sigma(\xivector)}^{-1}
\phi(\xivector)
\chi_E(L(\xivector))
\,d\lambda(\xivector)\,d\rhovector.
\end{equation}
Since all $\beta_n$ are assumed to be strictly $>\tfrac12$,
applying Cauchy-Schwarz to the integral with respect
to all $\rho_n$ gives an upper bound
\begin{equation}  \label{hexpression}
\lesssim
\int_{S\subset \Xi}
\prod_{n=1}^4 h_n(\xi_n)
\abr{\sigma(\xivector)}^{-1}
\phi(\xivector)
\chi_E(L(\xivector))
\,d\lambda(\xivector)
\end{equation}
where $h_n(\xi_n) = \norm{\tilde g_n(\xi_n,\cdot)}_{\lt(\reals^1)}
=\norm{g_n(\xi_n,\cdot)}_{\lt(\reals^1)}$.

According to Lemma~\ref{lemma:routine}, \eqref{hexpression} is
\begin{equation}
\lesssim 
|E|^{1/2} 
\prod_n\norm{g_n}_{\lt(\reals^2)}
\max_{S}\Big(\phi(\xivector)\abr{\sigma(\xivector)}^{-1} |\xivector|^{1/2}\Big),
\end{equation}
since the linear functional $L$ does not vanish identically on $\Xi$.

It remains to treat the region where $\max_n |\rho_n|\ge \tfrac18\abr{\sigma(\xivector)}$.
It is no loss of generality to restrict attention to the region where
$|\rho_4|=\max_n|\rho_n|$. An upper bound for the integral over this region is
\begin{multline} 
\int_{(\rho_1,\rho_2,\rho_3)\in\reals^3}\int_{S\subset \Xi}
\Big(\int_{\rho_4\in\reals}
\abr{\rho_1-\rho_2+\rho_3-\rho_4+\sigma(\xivector)}^{-1}
g_4(\xi_4,\rho_4)\,d\rho_4\Big)
\\
\prod_{n=1}^3 
\big(g_n(\xi_n,\rho_n)
\abr{\rho_n}^{-\beta_n}\big)
\,\abr{\sigma(\xivector)}^{-\beta_4}
\phi(\xivector)
\chi_E(L(\xivector))
\,d\lambda(\xivector)
\,\prod_{n=1}^3d\rho_n.
\end{multline}
Consider the contribution of a subregion in  which $\abr{\rho_4}$
is comparable to an arbitrary constant $\Lambda\gtrsim \abr{\sigma(\xivector)}$.
Then the innermost integral is the convolution
of $g_4(\xi_4,\cdot)$ with an $\lt(\reals)$ function,
evaluated at $\rho_1-\rho_2+\rho_3+\sigma(\xivector)$, so the entire
integral can be written as
\begin{multline} 
\int_{\reals^3}\int_{S\subset \Xi}
G_4(\xi_4, \rho_1-\rho_2+\rho_3+\sigma(\xivector))
\prod_{n=1}^3 
g_n(\xi_n,\rho_n)
\phi(\xivector) \chi_E(L(\xivector))
\,d\lambda(\xivector)
\,\prod_{n=1}^4 \abr{\rho_n}^{-\beta_n} d\rho_n
\end{multline}
where $\norm{G_4}_{\lt}\le C\norm{g_4}_{\lt}$.
$G_4(\xi_4, \rho_1-\rho_2+\rho_3+\sigma(\xivector))$
can be reexpressed as $\tilde G_4(\xi_4, \rho_1-\rho_2+\rho_3+
\xi_1^2-\xi_2^2+\xi_3^2)$ where $\norm{\tilde G_4}_{\lt}
=\norm{G_4}_{\lt}$.

This expression is identical to the expression \eqref{eq:subregion}
reached in the proof of Proposition~\ref{prop:strichartzA}, 
with the role played by $\abr{\rho_\nu}^{-\beta_\nu}$
in \eqref{eq:subregion}
now taken by $\abr{\rho_4}^{-\beta_4}$.
The rest of the analysis is identical to that of
\eqref{eq:subregion}, with the simplification
that here all $\beta_m$ are $>\tfrac12$.
\end{proof}

\section{Conclusion of the proof of Proposition~\ref{prop:PhifromX}}
\label{section:conclusionofPhifromX}

\begin{lemma}
Suppose that $s<0$, $r>-\tfrac14$, and $b>\tfrac12$. 
Let $\psi(\xivector) := \sum_{n=1}^4 (-1)^n \abr{\xi_n}^{2s}$.
Then
\begin{equation} 
\int_{\reals^4}
\int_\Xi
\prod_{n=1}^4 \Big(
g_n(\xi_n,\tau_n)
\abr{\xi_n}^{-r}\abr{\tau_n-\xi_n^2}^{-b}
\Big)
|\psi(\xivector)|
\,d\lambda(\xivector) \,d\tauvector
\lesssim 
\prod_{n=1}^4 \norm{g_n}_{\lt(\reals^2)}.
\end{equation}
\end{lemma}

\begin{proof}
It is no loss of generality to assume throughout the proof that
$\norm{g_n}_{\lt}=1$ for all $n$. 
We analyze the integral \eqref{needStrichartz} using 
Proposition~\ref{prop:strichartzB},
with $\phi\equiv\psi$.
Recall the symmetries discussed in the proof of Proposition~\ref{prop:strichartzC}.
These will be used to reduce the number of cases that must be 
discussed in the proof.

Let $N\ge 1$, and
consider the contribution to the integral made
by the subregion $S_N$ of integration in which all $\abr{\xi_j}$ are comparable to $N$.
Because of the symmetries listed above, we may restriction attention
to the region where $|\xi_1-\xi_2|\le|\xi_1-\xi_4|$.
Consider the subregion $S_{N,A}$ where $|\xi_1-\xi_2|\sim AN$,
where $A\lesssim 1$, and examine the bound given by 
Proposition~\ref{prop:strichartzB} with $L(\xivector) = \xi_1-\xi_2$
and $E=[-CAN,CAN]$. 
Since
$|\psi(\xivector)|\lesssim N^{2s-2}|\sigma(\xivector)|$, 
the maximum value over this subregion of
$|E|^{1/2} |\psi(\xivector)| \abr{\sigma(\xivector)}^{-1}|\xivector|^{1/2}$
is 
$\lesssim (AN)^{1/2} N^{2s-2} N^{1/2} = A^{1/2}N^{2s-1}$.
Similarly
\begin{align*}
|E|^{1/2} |\xi_1-\xi_4|^{-1/2}|\psi(\xivector)| \abr{\sigma(\xivector)}^{-\beta}
&\lesssim
|E|^{1/2} |\xi_1-\xi_4|^{-1/2}|\psi(\xivector)| \abr{\sigma(\xivector)}^{-1/2}
\\
&\lesssim 
|\xi_1-\xi_2|^{1/2}
N^{2s-2}|\sigma(\xivector)|^{1/2} 
|\xi_1-\xi_4|^{-1/2}
\\
&= N^{2s-2}|\xi_1-\xi_2|
\\
&\lesssim AN^{2s-1}
\end{align*}
for all $\xivector\in S_{N,A}$.

Summing over dyadic values of $A\lesssim 1$ gives a bound of $\lesssim N^{2s-1}$. 
Taking the factors $\abr{\xi_n}^{-r}$ into account yields a net bound 
of $\lesssim N^{2s-4r-1}$ for
the contribution of $S_N$ to \eqref{needStrichartz}. 
Provided that
$-r<\tfrac14 -\tfrac12 s$,
this is $\lesssim N^{-\delta}$ for some $\delta>0$ and hence
we can sum over dyadic values of $N\ge 1$
to majorize the contribution of the entire region
on which all four quantities $\abr{\xi_n}$ are mutually comparable.
Since $s<0$, this is a less stringent condition on $r$ than 
the hypothesis $-r<\tfrac14$.
\qed

The relation $\xi_1-\xi_2+\xi_3-\xi_4=0$ defining $\Xi$ implies that
the largest two of the four quantities $|\xi_n|$  must remain uniformly comparable.
Consider next the contribution of a region of integration 
in which some two variables $\xi_n$ of opposite parity
are large, and at least one of the other two variables is comparatively small.
Because of symmetries, it is then no loss of generality to restrict
attention to the region where
$|\xi_1|,|\xi_2|\sim N_2$,
$\abr{\xi_3}\sim N_1$,
and $\abr{\xi_4}\sim N_0$,
where the parameters $N_0,N_1,N_2\ge 1$  satisfy
$N_0\le N_1\le N_2$.

In the subcase in which $N_0\sim N_1$,
consider the subregion $S_\Delta$
where $|\xi_4-\xi_3|\equiv|\xi_1-\xi_2|$ has some fixed order 
of magnitude $\Delta$; necessarily
$\Delta\lesssim N_1$.
There
$|\sigma(\xivector)| = |\xi_1-\xi_4|\cdot|\xi_1-\xi_2| \sim \Delta N_2$,
and
\begin{equation}
|\psi(\xivector)|\le 
|\varphi(\xi_3)-\varphi(\xi_4)|
+ |\varphi(\xi_1)-\varphi(\xi_2)|
\lesssim N_0^{2s-1}\Delta + N_2^{2s-1}\Delta
\lesssim N_0^{2s-1}\Delta
\end{equation}
since $|\xi_1-\xi_2|=|\xi_3-\xi_4|$.

Apply Proposition~\ref{prop:strichartzB}
with $L(\xivector) = \xi_4-\xi_3$
and $\phi(\xivector)=|\psi(\xivector)|$
to the contribution made by the region of integration
$S_\Delta$ to \eqref{needStrichartz}.
The maximum value over all $\xivector\in S_\Delta$ of
$ |L(\xivector)|^{1/2} |\psi(\xivector)| \abr{\sigma(\xivector)}^{-1}|\xivector|^{1/2} $
is
\begin{equation}
\lesssim
\Delta^{1/2}\cdot N_0^{2s-1}\Delta \cdot \abr{\Delta N_2}^{-1} N_2^{1/2}
\lesssim 
\Delta^{1/2} N_0^{2s-1} N_2^{-1/2}
= (\Delta/N_0)^{1/2} N_0^{2s-\tfrac12}N_2^{-1/2}.
\end{equation}
The maximum value over all $\xivector\in S$ of
$ |L(\xivector)|^{1/2} |\xi_1-\xi_4|^{-1/2}|\psi(\xivector)| 
\abr{\sigma(\xivector)}^{-\beta}$
is
\begin{equation}
\lesssim
\Delta^{1/2}\cdot N_2^{-1/2} \cdot N_0^{2s-1}\Delta\cdot \abr{\Delta N_2}^{-1/2}
\lesssim (\Delta/N_0) N_0^{2s} N_2^{-1}.
\end{equation}
Since $\Delta\lesssim N_0\lesssim N_2$,
the maximum of these two maxima  is 
$\lesssim (\Delta/N_0)^{1/2}N_0^{2s-\tfrac12}N_2^{-1/2}$.
Incorporating the factors $\abr{\xi_n}^{-r}$ from \eqref{needStrichartz}
introduces an additional factor of $N_0^{-2r}N_2^{-2r}$,
leaving a net bound of
$\lesssim (\Delta/N_0)^{1/2}N_0^{2s-\tfrac12-2r}N_2^{-\tfrac12-2r}$.
Summing over dyadic values of $\Delta\lesssim N_0$ yields a bound of
$\lesssim N_0^{2s-\tfrac12-2r}N_2^{-\tfrac12-2r}$
for the original region.
This quantity is $\lesssim N_2^{-\delta}$ for some $\delta>0$
if (and only if) $-r<\tfrac14$. 
We may then sum over dyadic $N_0\lesssim N_2$,
then over all dyadic $N_2$.
\qed

If on the other hand $N_0\le \tfrac1{10}N_1$ then $\Delta\sim N_1$ 
and $|\psi(\xivector)|\lesssim N_0^{-2s}$, so 
the maximum value of
$|L(\xivector)|^{1/2} |\psi(\xivector)| \abr{\sigma(\xivector)}^{-1}|\xivector|^{1/2}$
is 
$\lesssim N_1^{1/2}\cdot N_0^{2s}\cdot \abr{N_1N_2}^{-1}N_2^{1/2} $,
giving a net bound of
$N_0^{2s-r} N_1^{-\tfrac12-r}N_2^{-\tfrac12-2r}$,
which again is $\lesssim N_2^{-\delta}$ for some $\delta>0$
if and only if $-r<\tfrac14$.
Likewise the maximum value of
$ |L(\xivector)|^{1/2} |\xi_1-\xi_4|^{-1/2}|\psi(\xivector)| 
\abr{\sigma(\xivector)}^{-\beta}$
is
\begin{equation}
\lesssim \
N_1^{1/2}\cdot N_2^{-1/2}\cdot N_0^{2s}\cdot (N_1N_2)^{-1/2}
\lesssim
N_0^{2s}N_1^{0}N_2^{-1},
\end{equation}
leading once again to the less stringent requirement $-r<\tfrac14-\tfrac12 s$.
\qed

Because the roles of the four variables $\xi_n$ are not completely symmetric,
it is necessary to analyze separately the subcase in which
again $N_0\le \tfrac1{10}N_1\le \tfrac1{10}N_2$, but 
$\abr{\xi_4}\sim N_0$, $\abr{\xi_2}\sim N_1$, and $|\xi_1|,|\xi_3|\sim N_2$.
Then $|\sigma(\xivector)|\sim N_2^2$, and $|\psi(\xivector)|\lesssim N_0^{2s}$.
Thus $|\sigma(\xivector)|$ is at least as large as in the above analysis.
Since it was raised to negative powers above,
this new situation is more favorable. Therefore 
the hypothesis $-r<\tfrac14$ again suffices.
\qed

When the various symmetries between the indices $\{1,2,3,4\}$
are taken into account,
the above discussion exhausts all possible cases, and the proof is complete.
\end{proof}

\begin{proof}[Proof of Proposition~\ref{prop:PhifromX}]
It suffices to bound $\norm{u(t,\cdot)}_{H^s}$
for $t$ in the support of $\zeta_0$,
since $\scripti(t)\equiv 0$ for other $t$.
For such $t$, 
$u(t,x)\equiv \zeta_1(t)u(t,x)$ 
and hence $\widehat{u}$ can be replaced by $\widehat{\zeta_1(t)u}$ throughout
the above discussion.
Thus $\norm{u}_{X^{r,b}}$ can be replaced by
$\norm{\zeta_1 u}_{X^{r,b}}$ 
on the right-hand side of the inequality.
\end{proof}

\section{$Y^{s,b}$ norms}  \label{section:Ysb}

The purpose of this section is to introduce certain function spaces
$Y^{s,b}$, variants of the spaces $X^{s,b}$ employed by Bourgain \cite{bourgain}
and then Kenig, Ponce, and Vega \cite{kpv1} to establish wellposedness
of the nonlinear Schr\"odinger and Korteweg-de Vries equations.
An {\em a priori} bound for $|u|^2 u$ in these spaces, in terms of $u$,
will be proved in the following section.

Proposition~\ref{prop:PhifromX} asserts an {\em a priori} upper bound 
for a solution in $C^0(H^s)$ in terms of an $X^{r,b}$ bound. 
Rather than establishing an $X^{r,b}$ bound directly, we will work with $Y^{s,b}$.
Whereas the usual argument 
establishing an {\em a priori} $X^{0,b}$ bound for a solution 
breaks down for $X^{s,b}$ for $s$ strictly negative,
it continues to apply for $Y^{s,b}$ when an upper bound in $C^0(H^s)$ is known. 
$Y^{s,b}$ strictly contains $X^{s,b}$, but embeds in $X^{r,b}$ for certain $r<s$;
see Lemma~\ref{lemma:convertYtoX}.

Define the scaling operator
\begin{equation}
T_\lambda u(t,x) := \lambda u(\lambda^2 t,\lambda x);
\end{equation}
$T_\lambda$ acts on distributions $u$ defined on $\reals^2$.
It maps any solution of the cubic nonlinear Schr\"odinger equation to another solution. 
We use the same notation for functions of $x$ alone:
$T_\lambda f(x) := \lambda f(\lambda x)$.

Define also the (rough) Littlewood-Paley projections
\begin{equation}
\widehat{P_{<N} u}(\xi,\tau) := 
\begin{cases} \widehat{u}(\xi,\tau) \ &\text{if } |\xi|\le N
\\
0  &\text{if } |\xi|>N.
\end{cases}
\end{equation}
We say that a function $f$ is {\em $M$-band-limited}
if $\widehat{f}(\xi,\tau)=0$ whenever $|\xi|>M$.

Fix an infinitely differentiable, compactly supported 
cutoff function $\eta\in C^\infty_0(\reals^1)$ satisfying $\eta(0)\ne 0$.
\begin{definition}[$Y^{s,b}$ norm]\label{defn:Ysb}
Let $s,b\in\reals$ with $s\in[-\tfrac12,0]$. 
For any tempered distribution $u$ defined on
$\reals^2$  whose space-time Fourier transform $\widehat{u}(\xi,\tau)$
belongs to $L^2_{\text{loc}}(\reals^2)$,
\begin{equation}
\Norm{u}_{Y^{s,b}}
:= \sup_{t_0\in\reals} 
\sup_{N\ge 1} 
\Norm{ \eta(t-t_0)T_{N^{2s}}(P_{<N} u) }_{X^{0,b}}.
\end{equation}
\end{definition}
It would be slightly more natural to form an $\ell^2$ norm
over a dyadic sequence of values of $N$,
rather than a supremum, but the definition used here is a bit
simpler to work with, and is sufficient for our purpose.
Observe that if $f$ is $N$-band-limited,
then $T_{N^{2s}}P_{<N} f$ is $N^{1+2s}$-band-limited.

For functions $f$ supported
in any fixed bounded interval with respect to time $t$, 
\begin{equation} \label{variantYsb}
\sup_N \int_{\langle\xi\rangle\sim N} \int_{\tau\in\reals} |\widehat{f}(\xi,\tau)|^2 
\,\abr{\xi}^{2s}
\,\abr{N^{4s}(\tau-\xi^2)}^{2b}
\,d\xi\,d\tau
\lesssim \norm{f}_{Y^{s,b}}^2,
\end{equation}
although the reverse inequality does not hold;\footnote{
For $r=(1+4b)s$ 
and $|\xi|$ of some fixed order of magnitude $N\ge 1$,
the left-hand side of \eqref{variantYsb} is equivalent 
to the $X^{r,b}$ norm squared
in the region where $|\tau-\xi^2|\gtrsim N^{-4s}$;
it becomes larger as $|\tau-\xi^2|$ becomes smaller than this threshold.}
this inequality can be derived as in the proof 
of Lemma~\ref{lemma:convertYtoX} below.
Because $s$ is negative and $b$ positive, the factor
$\abr{N^{4s}(\tau-\xi^2)}^{2b}$ is weaker than the corresponding factor
$\abr{\tau-\xi^2}^{2b}$ that appears in the $X^{s,b}$ norm.
Thus $X^{s,b}$ embeds continuously in $Y^{s,b}$.

Our first lemma is a simple consequence of the definition; the proof is omitted.
\begin{lemma}[Insensitivity to smooth cutoffs]
(i) If $h:\reals\to\complex$ is compactly supported and infinitely differentiable
then
$\norm{hu}_{Y^{s,b}}
\lesssim \norm{u}_{Y^{s,b}}$
for all $u\in Y^{s,b}$.

(ii) Changing the cutoff function $\eta$ in the definition of $Y^{s,b}$
leads to an equivalent norm,
provided that $\eta\in C^\infty$ is compactly supported, and not identically zero.
\end{lemma}

\begin{remark} 
For $s<0$,
the spaces $Y^{s,b}$ are natural from the point of view of the 
extant $H^0$ theory. 
If an initial datum $u_0$ for \eqref{nls} is $N$-band-limited in the
sense that $\widehat{u_0}(\xi)$ is supported where $|\xi|\sim N$,
and if $\Norm{u_0}_{H^s}\sim 1$, then $u_0\in H^0$,
but with large norm $\norm{u_0}_{H^0}\sim N^{-s}$.
Hence the Cauchy problem with initial datum $u_0$ has a solution
belonging to $X^{0,b}$. This does not follow from the usual fixed point argument, since
$u_0$ may be quite large in $H^0$. 
Instead one can partition the interval $[0,t]$ into sufficiently short 
subintervals that a fixed point argument applies on each, and 
invoke $H^0$ norm conservation.

An equivalent way to do the first time step is to
solve the Cauchy problem for unit time
with rescaled initial datum $T_{\lambda_N}u_0$,
where $\lambda_N = N^{2s}$, then to reverse the scaling. 
The exponent is chosen so that
$\|T_{\lambda_N}u_0\|_{H^0} \lesssim 1$ uniformly in
$N\ge 1$. 
Successive time steps are done in the same way. 
\end{remark}

The next simple lemma makes possible the conversion of bounds in $Y^{s,b}$ to the
more standard spaces $X^{r,b}$. 

\begin{lemma}[$Y$ controls $X$]\label{lemma:convertYtoX}
Let $s<0$ and $b\ge 0$.
For any $A<\infty$ and
any $r<(1+4b)s$ and all Schwartz class functions $f(t,x)$
supported where $|t|\le A$, we have
\begin{equation}
\Norm{f}_{X^{r,b}} \lesssim \Norm{f}_{Y^{s,b}}.
\end{equation}
\end{lemma}

The converse inequality is not true; 
in the region where 
$|\tau-\xi^2|\ll \abr{\xi}^{-4s}$,
the $Y^{s,b}$ norm is stronger than the $X^{r,b}$ norm
even for $r=(1+4b)s$.
We make this conversion both for the sake of conceptual simplicity,
and because it simplifies certain calculations later on;
retaining the full strength of the $Y^{s,b}$ bound does not
seem to lead directly to any improvement in our main theorems,
although it might contribute to some small improvement if combined
with other refinements.

While Lemma~\ref{lemma:convertYtoX} is needed to control $d\Phi/dt$, a 
variant will be used in establishing the $Y^{s,b}$ norm bound.
For any real number $M\ge 1$
define the $X^{r,b}_M$ and $\hat X^{r,b}_M$ norms by
\begin{align*}
\Norm{f}_{X^{r,b}_M}^2
&:=\iint_{\reals^2}
|\widehat{f}(\xi,\tau)|^2 \abr{\xi/M}^{2r} \abr{\tau-\xi^2}^{2b}\,d\tau\,d\xi.
\\
\Norm{g}_{\hat X^{r,b}_M}^2
&:=\iint_{\reals^2}
|{g}(\xi,\tau)|^2 \abr{\xi/M}^{2r} \abr{\tau-\xi^2}^{2b}\,d\tau\,d\xi
= c\Norm{\check{g}}_{X^{r,b}_M}^2.
\\
\intertext{Likewise define}
\Norm{g}_{\hat X^{r,b}}^2
&:=\iint_{\reals^2}
|{g}(\xi,\tau)|^2 \abr{\xi}^{2r} \abr{\tau-\xi^2}^{2b}\,d\tau\,d\xi.
\end{align*} 

\begin{lemma}[$Y$ controls $X$, again]\label{lemma:technicalembedding}
Let $s<0$ and $b\in(\tfrac12,1)$.
Let $r< (1+4b)s$. 
Then for any $f\in Y^{s,b}$, 
any $N\ge 1$,
and any $t_0\in\reals$,
the function
$g(t,x) = \eta(t-t_0) T_{N^{2s}} f(t,x)$ 
belongs to $X^{r,b}_{N^{1+2s}}$ with bound
\begin{equation}
\Norm{g}_{X^{r,b}_{N^{1+2s}}}
\lesssim
\Norm{f}_{Y^{s,b}}.
\end{equation}
\end{lemma}

Choose any smooth, compactly supported function $\eta$
such that $\sum_{j\in\integers} \eta(t-j)\equiv 1$ for all $t\in\reals$.

\begin{lemma}[Littlewood-Paley inequality]\label{lemma:almostorthog}
Let $s\le 0$ and $b\in\reals$.
Let $g$ be any Schwartz function, and define
$g_j = \eta(t-j)g$ so that $g=\sum_{j\in\integers} g_j$.
Then the summands $g_j$ are almost orthogonal in $X^{0,b}$ norm, in the sense that
\begin{equation} \label{almostorthogonality}
\bignorm{g}_{X^{0,b}}\le C\Big(\sum_j\Norm{g_j}_{X^{0,b}}^2 \Big)^{1/2}
\end{equation}
where $C<\infty$ depends only on $s,b,\eta$.
\end{lemma}

\begin{proof}
Introduce the spatial Fourier transform 
$\scriptf g(t,\xi) = \int_{\reals} g(t,x)e^{-ix\xi}\,dx$.
Let $J(t)$ be the distribution in $\scripts'(\reals^1)$
whose Fourier transform is $\abr{\tau}^b$. Then $J$ may be decomposed as
$J = J_0+J_\infty$
where $J_0$ is compactly supported and $J_\infty$ belongs
to the Schwartz class.

Now 
\begin{equation}
\Norm{g}_{X^{0,b}} = \Norm{\scriptf g*(e^{i\xi^2t}J(t))}_{L^2}
\end{equation}
where $*$ denotes convolution, taken with respect to the $t$ variable alone
for each fixed value of $\xi$.
Since $J_\infty$ is a Schwartz function,
\begin{equation}
\Norm{\scriptf g*(e^{i\xi^2t}J_\infty(t))}_{L^2}
\lesssim (\sum_j \Norm{g_j}_{L^2}^2 )^{1/2},
\end{equation}
and since $b\ge 0$,
$\Norm{g_j}_{L^2}
\lesssim \Norm{g_j}_{X^{0,b}}$.

There exists a finite constant $C_0$,
depending only on $\eta$ and on the support of $J_0$,
such that no point 
$(t,x)$ belongs to the support of $g_j$ for more than $C_0$ integers $j$.
Because the cutoff functions $\eta(t-j)$ are independent of $x$,
the same bounded overlap property holds for their spatial Fourier
transforms $\scriptf g_j(t,\xi)$.
Because $J_0$ has compact support, 
it follows that likewise no point
$(t,\xi)$ belongs to the support of
$\scriptf g_j*(e^{i\xi^2t}J_0(t))$
for more than $C_0$ integers $j$.  

Therefore
\begin{align*}
\Norm{\scriptf g*(e^{i\xi^2t}J_0(t))}_{L^2}^2
&\lesssim \sum_j \Norm{\scriptf g_j*(e^{i\xi^2t}J_0(t))}_{L^2}^2 
\\
&\lesssim \sum_j \iint |\widehat{g_j}(\tau,\xi)|^2 |\widehat{J_0}(\tau-\xi^2)|^2
\,d\tau\,d\xi
\\
&\lesssim \sum_j\norm{g_j}_{X^{0,b}}^2
\end{align*}
since $|\widehat{J_0}| = |\widehat{J}-\widehat{J_\infty}|
\le |\widehat{J}|+C\lesssim \abr{\tau}^b+C \le \abr{\tau}^b$
since $b\ge 0$.
\end{proof}

\begin{proof}[Proof of Lemma~\ref{lemma:convertYtoX}]
Let $f$ be given.
Let $r := (1+4b)s$. It suffices to show that for all $N\ge 1$,
\begin{equation}
\int_{\abr{\xi}\sim N} \int_{\tau\in\reals}
|\widehat{f}(\xi,\tau)|^2 \abr{\tau-\xi^2}^{2b} \abr{\xi}^{2r}
\,d\xi\,d\tau
\lesssim \Norm{f}^2_{Y^{s,b}},
\end{equation}
since summing over all $N= 1,2,4,8,\dots$ then yields the required bound
for all $r$ strictly less than $(1+4b)s$.

Define $g_j := \eta(t-j)\cdot T_{N^{2s}}P_{<N} f$,
and $g:=\sum_{j\in\integers} g_j$,
as in Lemma~\ref{lemma:almostorthog}.
All but at most $CN^{-4s}$
terms in this decomposition vanish identically, 
because of the hypothesis restricting
the support of $f$ with respect to $t$.
Moreover 
$\widehat{f}(\xi,\tau) = N^{4s}\widehat{g}(N^{2s}\xi,N^{4s}\tau)$. 
Consequently a trivial majorization of the $\ell^2$ outer norm in 
\eqref{almostorthogonality} gives
\begin{equation} \label{ellinfinitybound}
\Norm{g}_{X^{0,b}}\lesssim 
N^{-2s}\max_j\Norm{g_j}_{X^{0,b}}
\lesssim N^{-2s}\Norm{f}_{Y^{s,b}}.
\end{equation}

Now (since $1+2s>0$)
\begin{align*}
\int_{\abr{\xi}\sim N} \int_{\tau\in\reals}
&|\widehat{f}(\xi,\tau)|^2 \abr{\tau-\xi^2}^{2b} \abr{\xi}^{2r}
\,d\xi\,d\tau
\\
&=
N^{8s}
\int_{\abr{\xi}\sim N} \int_{\tau\in\reals}
|\widehat{g}(N^{2s}\xi,N^{4s}\tau)|^2 \abr{\tau-\xi^2}^{2b} \abr{\xi}^{2r}
\,d\xi\,d\tau
\\
&=
N^{2s}
\int_{\abr{\xi}\sim N^{1+2s}} \int_{\tau\in\reals}
|\widehat{g}(\xi,\tau)|^2 \abr{N^{-4s}(\tau-\xi^2)}^{2b} \abr{N^{-2s}\xi}^{2r}
\,d\xi\,d\tau
\\
&\sim
N^{2s}
\int_{\abr{\xi}\sim N^{1+2s}} \int_{\tau\in\reals}
|\widehat{g}(\xi,\tau)|^2 \abr{N^{-4s}(\tau-\xi^2)}^{2b} N^{2r}
\,d\xi\,d\tau
\\
&\lesssim 
N^{2s-8bs+2r} 
\int_{\abr{\xi}\sim N^{1+2s}} \int_{\tau\in\reals}
|\widehat{g}(\xi,\tau)|^2 \abr{\tau-\xi^2}^{2b} 
\,d\xi\,d\tau
\\
&\le
N^{2s-8bs+2r} \Norm{g}_{X^{0,b}}^2.
\\
&\lesssim N^{-2s-8bs+2r} \Norm{f}_{Y^{s,b}}^2
\end{align*}
by \eqref{ellinfinitybound}. 
This is $\lesssim \Norm{f}_{Y^{s,b}}^2$
under the hypothesis that $r\le (1+4b)s$.
\end{proof}

The proof of the related embedding property stated in
Lemma~\ref{lemma:technicalembedding}
is nearly identical to that of Lemma~\ref{lemma:convertYtoX},
so is left to the reader.
\qed

Proposition~\ref{prop:PhifromX} together with the embedding
of $Y^{s,b}$ in $X^{r,b}$ established in Lemma~\ref{lemma:convertYtoX} yield
\begin{proposition} \label{prop:C0Hsbound}
Let $T_0<\infty$, $T\in[0,T_0]$, $s\in (-\tfrac12,0)$, $b\in (\tfrac12,1)$.
For any sufficiently smooth solution $u$ of \eqref{nlsmodified}
with initial datum $u_0$,
\begin{equation}
\Norm{u}_{C^0([-2T,2T],H^s)}^2 \le \Norm{u_0}_{H^s}^2
+ C\Norm{\zeta_1 u}_{Y^{s,b}}^4
\end{equation}
provided that $s<0$, $b>\tfrac12$, and
$-s<\tfrac14(1+4b)^{-1}$.
\end{proposition}

To use this bound we of course need to control the $Y^{s,b}$ norm of $u$.  This will be accomplished
in the next two sections.

\section{Bound for $|u|^2 u$}

The objective of this section is to prove the following nonlinear estimate.

\begin{proposition}[Trilinear estimate in $Y^{s,b}$]\label{prop:multiplyYsb}
Suppose that $s>-\tfrac2{15}$ and $b\in(\tfrac12,1)$ satisfy
\begin{equation} \label{sbYrequirement}
-s<(1+4b)^{-1}
\min\Big( 
\tfrac1{10}+\tfrac35(1-b),
\tfrac1{12}+\tfrac23(1-b)
\Big).
\end{equation}
Then for any $u,v,w\in Y^{s,b}$, 
\begin{equation} \label{trilinearsbY}
\norm{u\bar v w}_{Y^{s,b-1}}\lesssim
\norm{u}_{Y^{s,b}}
\norm{v}_{Y^{s,b}}
\norm{w}_{Y^{s,b}}.
\end{equation}
\end{proposition}
The product $u\bar v w$,
by virtue of having a locally integrable space-time Fourier transform,
consequently has a natural interpretation as a distribution.

\eqref{trilinearsbY} is a variant of a well-known inequality 
in which $Y^{s,c}$ is replaced by $X^{0,c}$ throughout.
Here there is a tradeoff:
Once the parameter $N$ in the definition of $Y^{s,b-1}$ is fixed,
no bound is asserted for $\widehat{u\bar v w}(\xi,\tau)$
for $|\xi|\gg N$, but $u,v,w$ are allowed
to lie in spaces of mildly negative order.

The right-hand side of \eqref{sbYrequirement}
equals $\tfrac2{15}$ when $b=\tfrac12$.
Thus for any $s>-\tfrac2{15}$ there does exist
$b\in(\tfrac12,1)$ satisfying \eqref{sbYrequirement}.

\begin{proof}
The definition of the $Y^{s,b}$ norm involves a supremum over $N\ge 1$; fix $N$.
Set $M:=N^{1+2s}$.
Choose $r$ very slightly less than $(1+4b)s$,
and recall the $X^{r,b}_{N^{1+2s}}$ bound formulated in 
Lemma~\ref{lemma:technicalembedding}. 

Pair the space-time Fourier transform of $u\bar v w$
with $\abr{\tau-\xi^2}^{b-1}g_4(\xi,\tau)$
where $g_4\in \lt(\reals^2)$. Substitute for the Fourier
transforms of $u,v,w$ as in \eqref{substitutegsubn}. 
Matters then reduce to showing that
$$
\int_{\xivector\in\Xi}
\int_{\tauvector\in\Xi}
\prod_{n=1}^4 
\Big( g_n(\xi_n,\tau_n)
\abr{\xi_n/M}^{-r} \abr{\tau_n-\xi_n^2}^{-\beta_n}
\Big)
\chi_{S_0}(\xivector)
\,d\lambda(\tauvector)\,d\lambda(\xivector)
\lesssim \prod_{n=1}^4\norm{g_n}_{\lt(\reals^n)}
$$
uniformly for all $M\ge 1$, where
$\beta_n:=b$ for $n\le 3$ and $\beta_4:=1-b$,
and $S_0:=\{\xivector: |\xi_4|\lesssim M\}$.
Assume with no loss of generality
that $\norm{g_n}_{L^2(\reals^2)}=1$ for all indices $n$.

An important special case arises when all $|\xi_n|$ are $\lesssim M=N^{1+2s}$.
For this subregion, the desired inequality is nothing more than
the well-known $X^{0,b-1}$ bound for $|u|^2u$
in terms of $\norm{u}_{X^{0,b}}^3$ (see e.g. \cite{tao}).

Consider next the contribution to the integral of the region where
$|\xi_n|\sim AM$ for all $n\ne 4$ for some single $A\gg 1$. 
For all such $\xivector$, 
$|\sigma(\xivector)|\sim (AM)^2$, so 
since $\min(b,1-b)=1-b$,
an application of Proposition~\ref{prop:strichartzC}
with $L(\xivector)=\xi_4$ yields an upper bound of the form
\begin{equation}
\frac{M^{1/4}}{(AM)^{1/4}} (AM)^{-2(1-b)} A^{-3r}
= M^{-2(1-b)}A^{-\tfrac14 -2(1-b)-3r}
\end{equation}
and we need both exponents to be negative. The exponent $-2(1-b)$ on $M$ is
certainly negative since $b<1$. Thus we need 
\begin{equation} \label{firstYsbrequirement}
-r < \tfrac1{12} + \tfrac23 (1-b).
\end{equation}
$Y^{s,b}$ embeds in $X^{r,b}_M$ for all
$r<(1+4b)s$ uniformly in $M\ge 1$, in the sense expressed by 
Lemma~\ref{lemma:technicalembedding},
so this expression is appropriately controlled by the product of $Y^{s,b}$ 
norms provided that \eqref{sbYrequirement} is satisfied.

A more delicate case arises
when $|\xi_j|\sim AM$ with $A\gg 1$ for two values of $j\in\{1,2,3\}$,
but $|\xi_n|\sim BM$ where $B\le A/10$ for the third index.
If $n=2$, then $\sigma(\xivector)\sim (AM)^2$,
and the above analysis applies; the sole change is that one factor
of $A^{-r}$ is now merely $\lesssim B^{-r}$,
which is a more favorable bound since $B\le A$
and $r<0$.
Thus it remains only to discuss the case
where $n$ is odd; by virtue of the symmetries
of the problem, it is then no loss of generality to suppose that $n=3$.

In the subcase where $B\gtrsim 1$,
we have $|\sigma|\gtrsim AMBM$ and  
Proposition~\ref{prop:strichartzC}, again with $L(\xivector)=\xi_4$, 
yields the upper bound
\begin{equation}
\frac{
M^{1/4}}{(AM)^{1/4}} (ABM^2)^{-(1-b)}A^{-2r}B^{-r}
= M^{-2(1-b)}A^{-\tfrac14-2r-(1-b)}B^{-r-(1-b)}.
\end{equation}
Provided that $-r<1-b$, the exponent on $B$ is negative, so
when $B\gtrsim A^{1/2}$ this is
$\lesssim M^{-2(1-b)}A^{-\tfrac14-\tfrac52 r -\tfrac32(1-b)}$.
In the case $1\lesssim B\lesssim A^{1/2}$ we invoke instead
Proposition~\ref{prop:strichartzA} with $L=\xi_4-\xi_3$ to obtain
an upper bound
\begin{equation}
\begin{aligned}
\frac{(BM)^{1/2}}{(AM)^{1/2}}
 (ABM^2)^{-(1-b)} A^{-2r}B^{-r}
&= M^{-2(1-b)} B^{\tfrac12-(1-b)-r} A^{-\tfrac12-(1-b)-2r}
\\
&\lesssim
M^{-2(1-b)} A^{-\tfrac14-\tfrac32(1-b)-\tfrac52 r}
\end{aligned}
\end{equation}
since the exponent $\tfrac12-(1-b)-r$ is positive for $b>\tfrac12$ and $r<0$,
and $B\lesssim A^{1/2}$.
This is the same bound as obtained above for $B\gtrsim A^{1/2}$.
The exponent on $M$ is negative since $b<1$, while
the exponent on $A$ is negative if
\begin{equation}
\label{secondYsbrequirement}
-r<\tfrac1{10}+\tfrac35(1-b).
\end{equation}
Under those conditions, this bound is 
summable over dyadic values of $M,A,B$. 

$1-b>\tfrac12>
\min( \tfrac1{10}+\tfrac35(1-b), \tfrac1{12}+\tfrac23(1-b))$ 
for all $b\in(\tfrac12,1)$, so the condition
that $-r<1-b$ does not appear in the hypotheses of the Proposition. 

If $ABM^2\lesssim 1$ then we use the upper bound
$\lesssim 1$ for $\abr{\sigma}$ in place of $(ABM^2)^{-(1-b)}$,
and obtain the upper bound
\begin{equation}
(BM)^{1/2}(AM)^{-1/2}  A^{-2r}
= B^{1/2}A^{-\tfrac12 -2r} 
\lesssim (A^{-1}M^{-2})^{1/2}A^{-\tfrac12 -2r} 
= A^{-1-2r}M^{-1}.
\end{equation}
Both exponents are negative for all $-r<\tfrac12$,
so this is a less stringent requirement than \eqref{secondYsbrequirement}.

Choosing $r$ to be sufficiently close to $(1+4b)s$ reduces
all these restrictions  to the stated hypothesis on $s$.
\end{proof}

\section{{\it A priori} bound in $Y^{s,b}$} \label{section:Ysbaprioribound}

The next result is the second main inequality underlying our theorems.

\begin{proposition} \label{prop:Ysbbound}
For any $s>-\tfrac2{15}$ and $b\in(\tfrac12,1)$
satisfying 
\begin{equation}
-s<(1+4b)^{-1}
\min\Big( \tfrac1{10}+\tfrac35(1-b), \tfrac1{12}+\tfrac23(1-b) \Big)
\end{equation}
any sufficiently smooth solution $u$ of \eqref{nlsmodified}
with initial datum $u_0$ satisfies
\begin{equation}
\label{firstcoupled}
\Norm{u}_{Y^{s,b}}
\lesssim \Norm{u}_{C^0(H^s)}
+ \Norm{u}_{Y^{s,b}}^3
\end{equation}
where $\Norm{\cdot}_{C^0(H^s)}:=\Norm{\cdot}_{C^0([-2T,2T],H^s)}$.
\end{proposition}

\begin{proof}
Choose $r<(1+4b)s$ sufficiently close to $(1+4b)s$.
Let $N\ge 1$, let $\eta$ be a smooth, compactly
supported function, and let $t_0\in\reals$.
Recall that $u$ may be considered to be defined, and to satisfy the 
modified equation \eqref{nlsmodified}, for all $t\in\reals$.

Consider $w(t,x) := \eta(t-t_0) T_{N^{2s}}(u)$,
which satisfies the equation
\begin{equation} \label{nlslocalized}
iw_t+w_{xx} = \eta'(t-t_0)T_{N^{2s}}u 
+ \eta(t-t_0)\zeta_0(N^{4s}(t-t_0))|T_{N^{2s}}u|^2 T_{N^{2s}}u.
\end{equation}
It suffices to bound $\widehat{w}(\xi,\tau)$ in the region where
$|\tau-\xi^2|\ge 1$, for the contribution of the 
region $|\tau-\xi^2|\le 1$ to the $X^{0,b}$ norm of $w$ is majorized
by $\lesssim \Norm{w}_{L^2(dt\,dx)}$, hence by $\lesssim \Norm{w}_{C^0(H^0)}$
because as a function of $t$, 
$w(t,x)$ is supported in an interval of uniformly bounded length;
hence this contribution is majorized by $\lesssim \Norm{u}_{C^0(H^s)}$.

We may express $\widehat{w}(\xi,\tau)$
as a constant times $(\tau-\xi^2)^{-1}$ times the Fourier transform
of the right-hand side of \eqref{nlslocalized}.
The contribution of the first term on the right is then easily handled, for
$\Norm{\eta'(t-t_0)T_{N^{2s}}u}_{L^2(dt\,dx)}
\le C\Norm{T_{N^{2s}}u}_{C^0(H^0)}
\le C\Norm{u}_{C^0(H^s)}$. 
After dividing by $\abr{\tau-\xi^2}^{-1}$
we therefore have a quantity whose norm in $X^{0,1}$
is majorized by 
$\lesssim \Norm{u}_{C^0(H^s)}$. 

The function $\eta(t-t_0)\zeta_0(N^{4s}(t-t_0))$
may be expressed as $\tilde\eta^3(t-t_0)$ where $\tilde\eta\in C^\infty$
is real-valued, is supported in a bounded interval independent of $N$, 
and is bounded in any $C^k$ norm uniformly in $N$.
The second term on the right-hand side of \eqref{nlslocalized}
thus becomes $|\tilde\eta(t-t_0)T_{N^{2s}}u|^2 \tilde\eta(t-t_0)T_{N^{2s}}u$.

By Lemma~\ref{lemma:technicalembedding},
the norm of $\tilde\eta(t-t_0)T_{N^{2s}}u$ in $X^{r,b}_{N^{1+2s}}$
is $\lesssim \Norm{u}_{Y^{s,b}}$.
Proposition~\ref{prop:multiplyYsb} 
says that the $X^{0,b}$ norm
of the function whose Fourier transform is $(\tau-\xi^2)^{-1}$
times the characteristic function of the region $|\xi|\lesssim N^{1+2s}$
times the space-time Fourier transform of
$|\tilde\eta(t-t_0)T_{N^{2s}}u|^2 \tilde\eta(t-t_0)T_{N^{2s}}u$
is majorized by
$\lesssim \Norm{u}_{Y^{s,b}}^3$,
provided that $-2+2b\le 1-2b$.
\end{proof}

\begin{proof}[Proof of Theorem~\ref{thm:aprioribound}]
For any finite $T$ and $\delta'>0$, there exists $\delta>0$ such that the bounds
of Propositions~\ref{prop:Ysbbound} and \ref{prop:C0Hsbound} together
imply an {\em a priori} upper bound 
$\norm{u}_{C^0([0,T],H^s)}\le\delta'$
provided that $\norm{u_0}_{H^s}\le\delta$
and 
$\norm{u}_{C^0([0,T],H^s)}\le 2\delta'$.

To prove the theorem,
it suffices to show that given any $R<\infty$, there exists $\eps_0>0$
such that for any $u_0\in H^0$ satisfying $\norm{u_0}_{H^s}\le R$,
if $u$ denotes the solution of \eqref{nls} with initial datum $u_0$, 
then $T_{\eps_0}u$ satisfies an {\em a priori} $C^0([0,1],H^s)$ bound.
Because $s>-\tfrac12$, the equation is subcritical in $H^s$; 
there exists $\eps_0$ so that $\norm{\eps u_0(\eps x)}_{H^s}\le\delta$
whenever $\norm{u_0}_{H^s}\le R$ and $0<\eps\le\eps_0$.
We know that $u\in C^0(H^0)$, hence $u\in C^0(H^s)$.
For very small $\eps$, depending on $\norm{u_0}_{H^0}$,
we have 
$\norm{T_\eps u}_{C^0([0,1],H^s)}\le\delta'$.

Now a continuity argument can be applied.
If $\eps>0$ has the property that
$\norm{T_\eps u}_{C^0([0,1],H^s)}\le \delta'$,
then there exists $\eps'>\eps$ such that
$\norm{T_{\eps'} u}_{C^0([0,1],H^s)}\le 2\delta'$,
and provided that $\eps'\le \eps_0$ and $\eps_0$ is
chosen to be sufficiently small but depending only on $R$, 
this implies that 
$\norm{T_{\eps'} u}_{C^0([0,1],H^s)}\le \delta'$.
Standard reasoning shows that this must then hold for $\eps'=\eps_0$.
\end{proof}

\section{Existence of weak solutions} \label{section:weaksolnsexist}

We now prove a weakened variant of Theorem~\ref{thm:weaksolns}
on the existence of weak solutions, showing merely that
weak solutions exist in $L^\infty(H^s)\cap C^0(H^{s'})\cap Y^{s,b}$
for all $s'<s$. The last detail, existence in $C^0(H^s)$,
will be addressed in \S\ref{section:equicontinuity}.

\begin{lemma} \label{lemma:equicontinuity}
\label{lemma:cheapequicontinuity}
Let $s>-\tfrac1{12}$.
Let $u_0\in H^s$, $\eps>0$, and $M<\infty$ be given.
There exist $T'>0$ and $\delta>0$ such that for any initial datum $v_0\in H^0$
satisfying $\Norm{v_0-u_0}_{H^s}<\delta$,
the standard solution $v$ of \eqref{nls} with initial datum $v_0$ satisfies
\begin{equation}
\int_{|\xi|\le M} |\widehat{v}(t_1,\xi)-\widehat{v}(t_2,\xi)|^2\abr{\xi}^{2s}\,d\xi
<\eps \text{ for all }  t_1,t_2\in[0,T']
\text{ satisfying } |t_1-t_2|<\delta.
\end{equation}
\end{lemma}

\begin{proof}
Fix any $b>\tfrac12$.
For any $\eps'>0$ there exists $\delta'>0$
such that any $w\in X^{0,b}$ satisfies
$\Norm{w(t_1,\cdot)-w(t_2,\cdot)}_{L^2}\lesssim |t_1-t_2|^\gamma \Norm{w}_{X^{0,b}}$
for all $\gamma<b-\tfrac12$ whenever $|t_1-t_2|\le 1$, as follows
from a standard Cauchy-Schwarz calculation.
By rescaling we conclude that 
\begin{equation}
\Norm{P_{<M} v(t_1,\cdot)-P_{<M} v(t_2,\cdot)}_{H^s}\le C_M|t_1-t_2|^\gamma \Norm{v}_{Y^{s,b}}
\end{equation}
whenever $|t_1-t_2|\lesssim M^{4s}$.

We have already established an {\em a priori\/} upper bound
for $\Norm{v}_{Y^{s,b}}$ in terms of $\norm{v_0}_{H^s}$,
hence in terms of $\norm{u_0}_{H^s}$ so long as $\delta\le 1$. 
Consequently
\begin{equation} \label{prepareforrellich}
\int_{|\xi|\le M} |\widehat{v}(t_1,\xi)-\widehat{v}(t_2,\xi)|^2 \abr{\xi}^{2s}\,d\xi
\le C'_M \eps'^2 
\end{equation}
provided that $|t_1-t_2|<\delta' M^{4s}$. The claim follows.
\end{proof}

\begin{proof}[Proof of Theorem~\ref{thm:weaksolns}]
Let $s\in (-\tfrac1{12},0)$, and then let
$s'\in (-\tfrac1{12},s)$ be arbitrary.
Consider any initial datum $u_0\in H^s$.
Let $(v_{0,j})$ be any sequence of functions in $H^0(\reals)$ 
such that $v_{0,j}\to u_0$ in $H^s$ norm as $j\to\infty$.
Let $v^{(j)}\in X^{0,b}$ be the unique standard 
solution of the Cauchy problem \eqref{nls}
with initial datum $v_{0,j}$. 

There exist $b>\tfrac12$ and  $T$ 
such that the sequence $v^{(j)}$ is uniformly bounded
in $C^0((-2T,2T),H^s)\cap Y^{s,b}$  norm.
Moreover, the mappings 
$(-2T,2T)\owns t\mapsto v^{(j)}(t,\cdot)\in H^{s'}$
are equicontinuous, by virtue of Lemma~\ref{lemma:cheapequicontinuity}
and the inequality
\begin{equation} \label{ssprimeholder}
\int_{|\xi|\ge M} |\widehat{f}(\xi)|^2\abr{\xi}^{2s'}\,d\xi
\le
CM^{2s'-2s}
\norm{f}_{H^s}^2.
\end{equation}

For any large $N$, decompose $v^{(j)}$ as 
\begin{equation*}
v^{(j)}= 
v^{(j)}_{N;\, \text{high}}
+
v^{(j)}_{N;\, \text{low}}
\end{equation*}
where $\widehat{v^{(j)}_{N;\, \text{low}}}(t,\xi):=\widehat{v^{(j)}}(t,\xi)$ 
when $|\xi|\le N$ and $:=0$ otherwise.
The equicontinuity 
of the mapping $t\mapsto v^{(n)}(t,\cdot)\in H^{s'}$ 
implies precompactness of $\{v^{(j)}_{N;\, \text{low}}\}$ in $C^0_t(C^\infty_x)$
for $x$ in every bounded region, for every $N$. A diagonal argument produces
a subsequence, denoted again by $v^{(j)}$, 
such that for every $N$, $v^{(j)}_{N;\, \text{low}}$ converges in 
the $C^0(C^\infty)$ topology in every bounded region. 
Since $v^{(j)}$ is uniformly bounded in $C^0(H^s)$,
there exists a distribution $u\in\scriptd'$ such that
$v^{(j)}\to u$ in the topology of $\scriptd'$. 

Equicontinuity, the uniform upper bound on $v^{(n)}$ in $C^0(H^s)\cap Y^{s,b}$,
and \eqref{ssprimeholder} together
ensure (possibly after passage to the limit of some subsubsequence)
that
$u\in C^0(H^{s'})\cap L^\infty(H^s)\cap Y^{s,b}$. 
It follows likewise
that $u(0,\cdot)\equiv u_0(\cdot)$.
The proof that the limit of some subsequence actually belongs to $C^0(H^s)$
will be completed in \S\ref{section:equicontinuity}.

It remains  to show that $u$ is a weak solution of the equation.
To simplify notation, denote the nonlinearity by $\scriptn(v):= |v|^2v$.
It follows directly from the above convergence that
$\scriptn(v^{(j)}_{N;\, \text{low}})$ converges
to $\scriptn(u_{N;\, \text{low}})$ in $C^0(C^\infty_{\text{loc}})$ for every $N$.

For any $\eps>0$ there exists $N$ such that
\begin{equation}
\norm{\scriptn(v^{(j)})-\scriptn(v^{(j)}_{N;\, \text{low}})}_{Y^{s',b-1}} 
\le\eps \ \text{ for all $j$. }
\end{equation}
This follows from the basic trilinear estimate, Proposition~\ref{prop:multiplyYsb},
since 
$v^{(j)}_{N;\, \text{high}}$ is arbitrarily small in $Y^{s',b}$
provided $N$ is sufficiently large, while the low part 
is bounded uniformly in $N$.
Likewise
$\scriptn(u)-\scriptn(u_{N;\, \text{low}})$
is $\le \eps$ in $Y^{s',b-1}$ for all $j$. 

These conclusions together imply that $\scriptn(v^{(j)})\to\scriptn(u)$
in the topology of $\scriptd'$.
Since $v^{(j)}$ is a solution of \eqref{nls},
it follows that $u$ is likewise a solution.
\end{proof}

\section{Continuity in time} \label{section:equicontinuity}

Since weak limits cannot be taken directly in spaces $C^0(H^s)$,  
some additional argument is needed to ensure
that the weak limits constructed above do belong to these spaces.
In this section we bridge that gap by establishing
a certain limited equicontinuity with respect to time. 

Recall the expressions $\Phi_\varphi(t,u) = \int_{\reals}
|\widehat{u}(t,\xi)|^2\varphi(\xi)\,d\xi$.
Additional control on the solution $u$ can obtained by
analyzing $\Phi_\varphi(t,u)$ for weights $\varphi$ which are more general than $\abr{\xi}^{2s}$,
and are specifically adapted to the initial datum $u_0$ (cf. the ``frequency envelopes'' used for instance in \cite{benjamin}).
We have actually proved the following statement more general than that
announced earlier.

\begin{lemma} \label{lemma:generalweights}
Let $s >-\tfrac1{12}$ and $s'\in (s,0)$.
For any nonnegative $C^2$ weight function $\varphi$ satisfying 
\begin{equation}
\varphi(\xi)\le \abr{\xi}^{2s'},
\qquad
\varphi'(\xi)\le \abr{\xi}^{2s'-1},
\qquad
\varphi''(\xi)\le \abr{\xi}^{2s'-2},
\end{equation}
for any initial datum $u_0\in H^0$,
the standard solution $u(t,x)$ of \eqref{nls}  satisfies
\begin{equation} \label{generalweightsconclusion}
\Big| \Phi_\varphi(t,u) - \Phi_\varphi(0,u)\Big|
\le C\Norm{u}_{Y^{s,b}}^4.
\end{equation}
\end{lemma}

From this can be extracted a high-frequency continuity result.
\begin{lemma} \label{lemma:highfrequencybound}
Let $s>-\tfrac1{12}$.
Let $u_{0}\in H^s$ and $\eps>0$ be given.
There exist $\delta>0$ and $N<\infty$ such that
for all $v_0\in H^0$ satisfying $\Norm{v_0-u_0}_{H^s}<\delta$,
the standard solution $v$ of \eqref{nls} with initial datum $v_0$ satisfies
\begin{equation}
\int_{|\xi|\ge N} |\widehat{v}(t,\xi)|^2\abr{\xi}^{2s}\,d\xi<\eps
\end{equation}
for all $t\in[0,T]$. 
\end{lemma}
Here the timespan $T\in (0,\infty)$ is fixed,
and it is assumed that $\norm{u_0}_{H^s}$
is sufficiently small that the proof of Theorem~\ref{thm:aprioribound}
applies to all smooth solutions with initial data
satisfying $\norm{v_0-u_0}_{H^s}\le\delta_0$, where $\delta_0$
depends on $T$.

\begin{proof}
Fix any exponent
$s'\in (s,0)$. 
Let $\eps>0$ be given.
Choose $M<\infty$ so that 
$\int_{|\xi|\ge M} |\widehat{u_0}(\xi)|^2\abr{\xi}^{2s}\,d\xi<\eps^2$.
Then
there exist a large parameter $M'\ge M$ and a weight function $\varphi$
satisfying the three inequalities hypothesized in Lemma~\ref{lemma:generalweights},
with exponent $s'$,
such that 
\begin{alignat}{2}
&\eps^{-1}\abr{\xi}^{2s}\ge \varphi(\xi)\ge \abr{\xi}^{2s} &\qquad &\text{ for all $\xi$,}
\\
&\varphi(\xi)= \eps^{-1} \abr{\xi}^{2s} &&\text{ for all $|\xi|\ge M'$,}
\\
&\varphi(\xi)=\abr{\xi}^{2s} &&\text{ for all $|\xi|\le M$.}
\end{alignat}
$M',\varphi$ depend on $\eps$ and on  $s'$. The conclusion 
\eqref{generalweightsconclusion} of Lemma~\ref{lemma:generalweights} holds
with a constant $C$ independent of $M,\eps$.

Thus by \eqref{generalweightsconclusion},
\begin{equation}
\begin{aligned}
\int_{|\xi|\ge M'} |\widehat{v_0}(\xi)|^2\varphi(\xi)\,d\xi
&\le
2\int_{|\xi|\ge M'} |\widehat{u_0}(\xi)|^2\varphi(\xi)\,d\xi
+
2\int_{|\xi|\le M'} |\widehat{v_0}(\xi)-\widehat{u_0}(\xi)|^2\varphi(\xi)\,d\xi
\\
&\le 
2\eps
+ C_\varphi \norm{v_0-u_0}_{H^s}^2
\end{aligned}
\end{equation}
where $C_\varphi$ depends on $\varphi$, hence ultimately on $\eps$.
Therefore there exists $\delta>0$ such that
\begin{equation}
\int_{\reals} |\widehat{v_0}(\xi)|^2\varphi(\xi)\,d\xi
\le 3\eps
\end{equation}
for every $v_0\in H^0$ satisfying $\Norm{v_0-u_0}_{H^s}<\delta$.

For such initial data $v_0$,
the associated solutions $v$ have uniformly bounded $Y^{s,b}$ norms,
with a bound independent of $\eps$,  provided that $\delta$ is sufficiently
small.  Therefore by Lemma~\ref{lemma:generalweights},
$\Phi_\varphi(t,v)= 
\int_\reals |\widehat{v}(t,\xi)|^2\varphi(\xi)\,d\xi$
is bounded by a finite constant independent of $\eps,M,M'$
uniformly for all $t\in[0,T]$.
Therefore since $\abr{\xi}^{2s}\le\eps\varphi(\xi)$ for all $|\xi|\ge M'$,
\begin{equation}
\int_{|\xi|\ge M'} |\widehat{v}(t,\xi)|^2\abr{\xi}^{2s}\,d\xi
\le 
\eps \int_{|\xi|\ge M'} |\widehat{v}(t,\xi)|^2\varphi(\xi)\,d\xi
\lesssim \eps,
\end{equation}
provided that $t\in[0,T]$ and $\Norm{v_0-u_0}_{H^s}<\delta$.
\end{proof}

Thus if $u_0,v^{(j)}_0$ are initial data
with $u_0\in H^s$ and $v^{(j)}_0\in H^0$,
and if $v^{(j)}\to u_0$ in the $H^s$ norm, then
the corresponding standard solutions
$v^{(j)}$ 
form an equicontinuous family in $C^0(H^s)$.
Therefore passage to the limit through an appropriate subsequence
produces a solution in $C^0(H^s)$, satisfying the other conclusions
of Theorem~\ref{thm:weaksolns}.

\begin{remark}
Lemma~\ref{lemma:highfrequencybound}
has the following direct consequence.
Let $s>-\tfrac1{12}$. If there exists $r>-\infty$
for which the solution mapping from datum to solution of \eqref{nls}
is continuous from $H^s$ to $C^0([0,T],H^r)$,
then the solution mapping 
is continuous from $H^s$ to $C^0([0,T],H^s)$.
\end{remark}

\end{document}